\documentclass[11pt,notitlepage,twoside]{article}
\usepackage{latexsym}
\usepackage{amsmath}
\usepackage{amsfonts}
\usepackage{amssymb}
\renewcommand{\a }{\alpha }
\renewcommand{\b }{\beta }
\renewcommand{\d}{\delta }

\newcommand{\D }{\Delta }

\newcommand{\e }{\varepsilon }
\newcommand{\g }{\gamma}

\renewcommand{\l }{\lambda }

\newcommand{\n }{\nabla }
\newcommand{\var }{\varphi }

\newcommand{\s }{\sigma }

\renewcommand{\th }{\theta }

\renewcommand{\O }{\Omega }

\newcommand{\ov}{\overline}

\newcommand{\be}{\begin{equation}}
\newcommand{\ee}{\end{equation}}
\newenvironment{pf}{\noindent{\bf Proof.}\enspace}{
\hfill$\Box$\medskip}
\newenvironment{pfn}[1]{\noindent{\bf Proof of {#1}\enspace}}{
\hfill$\Box$\medskip}
\newcommand{\R}{\mathbb{R}}
 
\newtheorem{thm}{Theorem}[section]
\newtheorem{pro}[thm]{Proposition}
\newtheorem{lem}[thm]{Lemma}

\numberwithin{equation}{section}

\author{Khalil El Mehdi$^{a}$ \,\, and Mokhless Hammami$^{c}$\\
{\footnotesize
a:  Facult\'e des Sciences et Techniques, Universit\'e de Nouakchott, BP 5026, Nouakchott,}\\
{\footnotesize
 Mauritania and  The Abdus Salam ICTP, Trieste, Italy. E-mail : \texttt{khalil@univ-nkc.mr}}\\
{\footnotesize
b: D{\'e}partement de Math{\'e}matiques, Facult{\'e} des Sciences de Sfax, Route
Soukra, Sfax,}\\
{\footnotesize
Tunisia.  E-mail : \texttt{Mokhless.Hammami@fss.rnu.tn}}
}

\title { {\Large \textbf{Blowing up Solutions for a Biharmonic Equation  with
      Critical Nonlinearity } }}

\begin{document}

\date{ }

\maketitle
{\footnotesize

\noindent
{\bf Abstract.}
In this paper we consider the following  biharmonic equation with critical exponent  $(P_{\e})$ : $\D^2u=K u^{\frac{n+4}{n-4}-\e}$, $u>0$ in $\O$ and $u=\D u=0$ on $\partial\O$, where $\O$ is a smooth bounded domain in $\R^n$, $n\geq 5$, $\e$ is a small positive parameter, and $K$ is a smooth positive function in $\overline\O$. We construct solutions of $(P_{\e})$ which blow up and concentrate at strict local maximum of $K$ either at the boundary or in the interior of $\O$. We also construct solutions of  $(P_{\e})$  concentrating at an interior strict local minimum point of $K$. Finally, we prove a nonexistence result for the correponding supercritical problem which is in sharp contrast to what happened for $(P_\e)$.

}

\section{Introduction and   Results }
\mbox{}
In this paper, we are concerned with the concentration phenomena of the following biharmonic equation under the Navier boundary condition
$$
(P_{\e} ) \quad \left\{
\begin{array}{cc}
 \D ^2 u = K u^{p-\e},\,\, u>0 &\mbox{ in }\, \O \\
    \D u= u =0   \quad\quad     & \mbox{ on }\,  \partial  \Omega ,
\end{array}
\right.
$$
where $\O$ is a smooth bounded domain in $\R^n$, $n\geq 5$, $\e$ is a small positive parameter, $p+1= 2n/(n-4)$ is the critical Sobolev exponent of the embedding  $H^2(\O )\cap H^1_0(\O ) \hookrightarrow L^{p+1}(\O )$, and $K$ is a smooth positive function in $\overline\O$.

The study of concentration phenomena for second order elliptic equations involving nearly critical exponent has attracted considerable attention in the last decades. See for example  \cite{AP}, \cite{BLR}, \cite{BEGR}, \cite{BP}, \cite{CY}, \cite{DFM1}, \cite{DFM2}, \cite{H}, \cite{KR}, \cite{MiP}, \cite{MP}, \cite{R1}, \cite{R2}, \cite{R4} and the references therein.\\
However, as far as the authors know, the concentration phenomena for problem
$(P_\e)$ have been studied only in \cite{CG}, \cite{G} and \cite{BE} for $K\equiv 1$ only.

The purpose of the present paper is to construct solutions for $(P_\e)$
concentrating at various point of $\O$. More precisely, we are interested in
constructing solutions concentrating at a strict local maximum point of $K$
either at the boundary or in the interior of $\O$. We will also  construct
solutions concentrating at an interior strict local minimum point of
$K$. Similar results for Laplacian equation involving nearly critical Sobolev
exponent has been  proved by Chabrowski and Yan
\cite{CY}. Compared with the second order case, further technical difficulties
have to be solved by means of delicate and careful estimates. Our method 
 uses  some  techniques developed by  Bahri \cite{B},  Rey\cite{R1} and Ben Ayed-El Mehdi \cite{BE} in the framework of {\it Theory of critical points at infinity}. The main idea consists
in 
performing refined expansions of the Euler functional associated to our
variational problem, and
its gradient in a neighborhood of  potential concentration sets. Such
expansions are made possible through a finite dimension reduction argument.

To state our results, we need to introduce some notation.
We denote by $G$  the Green's function of $\D^2$, that is,
$$
\forall x\in\O \quad \left\{
\begin{array}{cc}
 \D ^2 G(x,.) = c_n\d_x &\mbox{ in }\, \O \\
    \D G(x,.)= G(x,.) =0      & \mbox{ on }\,  \partial  \Omega ,
\end{array}
\right.
$$
where $\d_x$ denotes the Dirac mass at $x$ and $c_n=(n-4)(n-2)|S^{n-1}|$.
We also denote by $H$ the regular part of $G$, that is,
$$
 H(x,y)=|x-y|^{4-n} -G(x,y),\,\,\mbox{for }\,\, (x,y)\in\O\times\O.
$$
Let
\begin{eqnarray}\label{e:11}
 \d _ {x,\l }(y) =  \frac {c_0\lambda
 ^{\frac{n-4}{2}}}{(1+\lambda^2|y-x|^2)^{\frac{n-4}{2}}},\,\,
 c_0=[(n-4)(n-2)n(n+2)]^{(n-4)/8},\,\,
 \l >0, \,\, x \in \R^n.
\end{eqnarray}
It is well known (see \cite{Lin}) that $\d_{x,\l}$ are the only
solutions of
\begin{eqnarray*}
 \D^2 u =  u^{\frac{n+4}{n-4}},\quad  u>0 \mbox{  in  } \R^n, \quad
  \mbox{with } u\in L^{p+1}(\R^n) \quad \mbox{and } \D u \in L^2(\R^n)
\end{eqnarray*}
\noindent
and are also the only minimizers of the Sobolev inequality on the
whole space, that is
\begin{eqnarray}\label{e:12}
 S =\inf\{||\D u||^{2}_{L^2(\R^n)}||u||^{-2}_{L^{\frac{2n}{n-4}}(\R^n)}
,\, s.t.\, \D u\in L^2 ,u\in L^{\frac{2n}{n-4}} ,u\neq 0 \}.
\end{eqnarray}
\noindent
 We denote by  $P\d _{x,\l}$ the projection of the $\d
_{x,\l}$'s onto $H^2(\O )\cap H^1_0(\O)$, defined by
$$ \D^2 P\d_{x,\l}=\D^2\d_{x,\l} \mbox{ in } \O  \mbox{ and }
\D P\d_{x,\l}=P\d_{x,\l}=0\mbox{ on } \partial \O
$$
and we set
$$
\var_{x,\l} = \d_{x,\l}- P\d_{x,\l}.
$$
The space $\mathcal{H}(\O) := H^2(\O )\cap H^1_0(\O)$ is equipped with the norm $||.||$ and its corresponding inner product $(.,.)$ defined by
\begin{align}
||u||&=\left(\int_\O |\D u|^2\right)^{1/2},\qquad u\in \mathcal{H}(\O ) \label{e:13}\\
(u,v)&=\int_\O \D u\D v,\qquad u,v \in \mathcal{H}(\O ). \label{e:14}
\end{align}
Let
\begin{align}
|u|_q&=|u|_{L^q(\O )} \label{e:15}\\
E_{x,\l}&=\{ v\in \mathcal{H}(\O)/(v,P\d_{x,\l})=(v,\frac{\partial
P\d_{x,\l}}{\partial\l})=(v,\frac{\partial P\d_{x,\l}}{\partial
x_j})=0,\,j=1,...,n\}\label{e:16}.
\end{align}

Now we state the main results of this paper.

\begin{thm}\label{t:11}
Let $x_0 \in \partial\O$ be a strict local maximum point of $K$ satisfying
\begin{eqnarray}\label{e:17}
K(x) \leq K(x_0) - a |x-x_0|^{2+\a} \quad\forall x \in B_\mu(x_0) \cap \overline\O,
\end{eqnarray}
where $\mu > 0$, $a>0$ and $\a\geq 0$ if $n\leq 6$, $\a \in [0, 4/(n-6))$ if $n\geq 7$. Then there is an $\e_0>0$, such that for each $\e\in (0,\e_0]$, $(P_\e)$ has a solution of the form
\begin{eqnarray}\label{e:18}
u_\e = \a_\e P\d _{x_\e , \l _\e}+v_\e,
\end{eqnarray}
where $v_\e \in E_{x_\e, \l_\e}$, and as $\e\to 0$,
\be\label{e:19}
\a_\e \to K(x_0)^{(4-n)/8},\,\,
||v_\e||\to 0,\,\,
x_\e \to x_0,\,\,
\l_\e \to + \infty,\,\,
\l_\e d(x_\e,\partial\O) \to +\infty.
\ee
\end{thm}
\begin{thm}\label{t:12}
Let $x_0 \in \O$ be a strict local maximum point of $K$.  Then there is an $\e_0>0$, such that for each $\e\in (0,\e_0]$, $(P_\e)$ has a solution of the form \eqref{e:18} satisfying \eqref{e:19}.
\end{thm}

The aim of the next result is to  show that  if $K$ is flat enough around a strict local minimum, $(P_\e)$ has a solution concentrating at this point.
\begin{thm}\label{t:14} Let $x_0\in\O$ be a strict local minimum point of $K$ satisfying
\begin{align}
&|D^l K(x)| \leq C |x-x_0|^{L-l},\quad l=1,...,n-4,\quad \forall x\in B_\mu (x_0), \label{e:114}\\
&|K(x) -K(x_0)|\geq C_0 |x-x_0|^L,\quad \forall x\in B_\mu (x_0), \label{e:115}
\end{align}
where $L>n-4$ is a constant, and where $C$ and $C_0$ are positive constants.\\
 Then there is an $\e_0>0$ such that for each $\e\in(0,\e_0]$, $(P_\e)$ has a solution of the form \eqref{e:18} satisfying \eqref{e:19} and $\e \l_\e^{n-4} \to c>0$.
\end{thm}

In the case $n=5$ or $6$, we can obtain a better result.
\begin{thm}\label{t:15}
Assume that $x_0 \in \O$ is a strict local minimum point of $K$. If one of the following conditions is satisfied :\\
$(i)$\quad $n=5$;\\
$(ii)$\quad $n=6$ and
\begin{eqnarray}\label{e:116}
c_1H(x_0,x_0) - \frac{c_2 \D K(x_0)} {36K(x_0)}> 0,\,\,
\mbox{with}\,\, c_1=c_0^{\frac{2n}{n-4}}\int_{\R^n}\frac{dy}{(1+
|y|^2)^{\frac{n+4}{2}}}\, c_2= \int_{\R^n} |y|^2\d_{o,1}^{p+1}dy,
\end{eqnarray}
then the conclusion of Theorem \ref{t:14} holds.
\end{thm}

The condition  \eqref{e:116} is nearly necessary. Indeed, we have the following result:
\begin{thm}\label{t:16}
Assume that $x_0 \in\O$ is a critical point of $K$ satisfying one of the following conditions :\\
$(i)$\quad $n\geq 7$ and $\D K(x_0) >0$,\\
$(ii)$ \quad $n=6$ and $c_1H(x_0,x_0) - \frac{c_2 \D K(x_0)} {36K(x_0)} < 0$, where $c_1$ and $c_2$
are the constants defined in Theorem \ref{t:15}.\\
Then $(P_\e)$ has no solution of the form \eqref{e:18} satisfying \eqref{e:19}.
\end{thm}

In contrast with the above results, we have the following nonexistence result for the supercritical problem.
\begin{thm}\label{t:17}
 Assume that $x_0 \in\O$ is a critical point of $K$ satisfying one of the following conditions :\\
$(i)$\quad $n=5$ ,\\
$(ii)$\quad $n=6$ and $c_1H(x_0,x_0) - \frac{c_2 \D K(x_0)}
{36K(x_0)} > 0$, where $c_1$ and $c_2$
 are the constants defined in Theorem \ref{t:15},\\
$(iii)$\quad $n\geq 7$ and $-\D K(x_0)> 0$.\\
 Then the problem 
 $$
 (Q_{\e} ) \quad \left\{
 \begin{array}{cc}
  \D ^2 u = K u^{p+\e},\,\, u>0 &\mbox{ in }\, \O \\
     \D u= u =0   \quad\quad     & \mbox{ on }\,  \partial  \Omega ,
 \end{array}
 \right.
 $$
 has no solution of the form \eqref{e:18} satisfying \eqref{e:19}.
 \end{thm}
The proof of our results is inspired by the
methods of  \cite{B}, \cite{BE}, \cite{CY} and \cite{R1}. The next section
will be devoted to some useful estimates needed in the proofs of our results. In section 3 we prove Theorems \ref{t:11}, \ref{t:12} and \ref{t:16}, while Theorems
\ref{t:14}, \ref{t:15} and \ref{t:17}   are proved in section 4. Lastly, we give in the appendix some integral estimates which are needed in Section 2.

\section{The Technical Framework}
\mbox{}
First of all, let us introduce the general setting.
For $\e>0$, we define on $\mathcal{H}\diagdown \{0\}$ the functional
\begin{eqnarray}\label{e:21}
J_\e (u) = \frac{\int_\O |\D u|^2}{\left(\int_\O K(x) |u|^{p+1-\e}\right)^{\frac{2}{p+1-\e}}}.
\end{eqnarray}
If $u$ is a critical point of $J_\e$, $u$ satisfies on $\O$ the equation
\begin{eqnarray}\label{e:22}
\D^2 u = l_\e(u) K(x) |u|^{p-1-\e}u
\end{eqnarray}
with
\begin{eqnarray}\label{e:23}
l_\e (u) = \frac{\int_\O |\D u|^2}{\int_\O K(x) |u|^{p+1-\e}}.
\end{eqnarray}
Conversely, we see that any solution of \eqref{e:22} is a critical point of $J_\e$.\\
Note that if $u$ is a positive critical point of $J_\e$, then $\left(l_\e(u)\right)^{\frac{1}{p-1-\e}} u$ is a solution of $(P_\e)$. This will allow us to look for solutions of $(P_\e)$ as critical points of $J_\e$.\\
Now let
$$
\mathcal{M}_\e = \{ (x,\l,v)\in \O \times \R^*_+ \times \mathcal{H}(\O)/ v \in E_{x,\l},\,\, ||v|| \leq \nu_0 \},
$$
where $\nu_0$ is a small positive constant.\\
Let us define the functional
\begin{equation}\label{psi}
\psi_\e : \mathcal{M}_\e \to \R, \,\, \psi_\e (x,\l,v) = J_\e
(P\d_{x,\l} +v).
\end{equation}
Notice that $(x,\l,v)$ is a critical point of $\psi_\e$ if and only if $u= P\d_{x,\l} + v$ is a critical point of $J_\e$. So this fact allows us to look for critical points of $J_\e$ by successive optimizations with respect to the different parameters on $\mathcal{M}_\e$.\\
On the other hand, $(x,\l,v) \in \mathcal{M}_\e$ is a critical point of $\psi_\e$ on $\mathcal{M}_\e$ if and only if there are $A$, $B$, $C_j \in \R$, $1\leq j\leq n$, such that
\begin{align*}
(E_{x_i}):& \quad \frac{\partial\psi_\e}{\partial x_i} = B\left( \frac{\partial ^2P\d_{x,\l}}{\partial\l\partial x_i}, v\right) + \sum_{j=1}^n C_j \left( \frac{\partial ^2P\d_{x,\l}}{\partial x_j\partial x_i}, v\right),\,i=1,...,n\\
(E_{\l}):& \quad \frac{\partial\psi_\e}{\partial \l} = B\left( \frac{\partial ^2P\d_{x,\l}}{\partial\l^2}, v\right) + \sum_{j=1}^n C_j \left( \frac{\partial ^2P\d_{x,\l}}{\partial x_j\partial \l}, v\right), \\
(E_{v}):& \quad \frac{\partial\psi_\e}{\partial v} = A P\d_{x,\l} +  B\frac{\partial P\d_{x,\l}}{\partial\l} + \sum_{j=1}^n C_j  \frac{\partial P\d_{x,\l}}{\partial x_j}.
\end{align*}
As usual in these types of problems, we first deal with the $v$-part of $u$. Namely, we prove the following.
\begin{pro}\label{p:21}
There exist $\e_1>0$, $\nu_0>0$, and a smooth map which to any $(\e, x, \l)\in (0,\e_1) \times \O \times \R^*_+$ with $\l d(x,\partial\O) > \nu_0^{-1}$, and $\e \log \l < \nu_0$, associates $v_\e= v_{\e,x,\l}\in E_{x,\l}$, $||v_\e|| < \nu_0$ such that $(E_v)$ is satisfied for some $(A,B,C_1,...,C_n)_{\e,x,\l}\in \R^{n+2}$. Such a $v_\e$ is unique, minimizes $\psi_\e(x,\l,v)$ with respect to $v$ in $\{v\in E_{x,\l}/ ||v||< \nu_0\}$, and we have the following estimate
$$
||v_\e||\ = O\left( \sum_{j=1}^k \frac{|D^jK(x)|}{\l^j} + \frac{1}{\l^{k+1}} + \e  + \frac{1}{(\l d)^{\frac{n-4}{2}+\th}}\right),
$$
where $\th >0$, $k$ is the biggest positive integer satisfying $k\leq \frac{n-4}{2}$, and where $d=d(x,\partial\O)$.
\end{pro}
\begin{pf}
As in \cite{B} (see also \cite{BC} and  \cite{R1}) we write
\begin{align}\label{e:24}
\psi_\e (x,\l,v)&=J_\e(P\d_{x,\l}+ v)\notag\\
&= \frac{||P\d_{x,\l} + v||^2}{\left(\int_\O K(y) |P\d_{x,\l}+v|^{p+1-\e}\right)^{2/(p+1-\e)}}\notag\\
&=\psi_\e(x,\l,0) -(f_\e,v) + \frac{1}{2}Q_\e(v,v) + O\left(||v||^{\min(3,p+1-\e)}\right),
\end{align}
where
$$
(f_\e,v)= 2J_\e(P\d_{x,\l})\frac{\int_\O K(y)P\d_{x,\l}^{p-\e}v}{\int_\O K(y)P\d_{x,\l}^{p+1-\e}},
$$
and
\begin{align*}
Q_\e(v,v)&= 2J_\e(P\d_{x,\l})\left[ \frac{||v||^2}{||P\d_{x,\l}||^2} - (p-\e)\frac{\int_\O K(y)P\d_{x,\l}^{p-1-\e}v^2}{\int_\O K(y)P\d_{x,\l}^{p+1-\e}}\right.\\
& \left.+ (p+3-\e)\left(\frac{\int_\O K(y)P\d_{x,\l}^{p-\e}v}{\int_\O K(y)P\d_{x,\l}^{p+1-\e}}\right)^2\right].
\end{align*}
It follows from Proposition 2.1 \cite{BH}, and Lemmas \ref{l:22} and \ref{l:23} that
\begin{align}
\int_\O K(y)P\d_{x,\l}^{p+1-\e} &= K(x)S_n + O\left(\frac{1}{\l^2}+ \e\log\l+ \frac{1}{(\l d)^{n-4}}\right),\label{e:25}\\
\int_\O K(y) P\d_{x,\l}^{p-\e}v& = O\left(\e + \sum_{j=1}^{k}
\frac{|D^jK(x)|}{\l^{j}} + \frac{1}{\l^{k+1}} + \frac{1}{(\l
d)^{\frac{n-4}{2}+ \th}}\right)||v||,\label{e:26}
\end{align}
where $k$ denotes the biggest positive integer satisfying $k\leq n-4/2$ and $\th >0$.\\
Now, we observe that
\begin{align}
\int_\O K(y)P\d_{x,\l}^{p-1-\e}v^2 &=\int_\O K(y)\d_{x,\l}^{p-1-\e}v^2 + o(||v||^2)\notag\\
&=\int_\O K(y)\d_{x,\l}^{p-1}v^2 + o(||v||^2)\notag\\
&=K(x)\int_\O \d_{x,\l}^{p-1}v^2 + o(||v||^2).\label{e:27}
\end{align}
One can check that (see \cite{BH})
\begin{eqnarray}\label{e:28}
||P\d_{x,\l}||^2 = S_n + O((\l d)^{4-n}).
\end{eqnarray}
Combining \eqref{e:25},..., \eqref{e:28}, we obtain
\begin{eqnarray}\label{e:29}
Q_\e(v,v)= \frac{2J_\e(P\d_{x,\l})}{S_n} \left[ ||v||^2 - p \int_\O \d_{x,\l}^{p-1}v^2 + o(||v||^2)\right].
\end{eqnarray}
According to \cite{BE1}, there exists some positive constant independent of $\e$, for $\e$ small enough, such that
\begin{eqnarray}\label{e:210}
 ||v||^2 - p \int_\O \d_{x,\l}^{p-1}v^2 \geq c ||v||^2,\quad \forall v\in E_{x,\l}.
\end{eqnarray}
It follows from Lemma \ref{l:23} that
\begin{eqnarray}\label{e:211}
(f_\e,v)= O\left(\sum_{j=1}^k \frac{|D^j K(x)|}{\l^j} +
\frac{1}{\l^{k+1}}+ \e + \frac{1}{(\l
d)^{\frac{n-4}{2}+\th}}\right).
\end{eqnarray}
It is easy to see that Proposition \ref{p:21} follows from \eqref{e:24},..., \eqref{e:211}.
\end{pf}
Next, we prove a useful expansion of the functional $J_\e$ associated to $(P_\e)$, and its gradient in a neighborhood of potential concentration sets.
\begin{pro}\label{p:24}
Suppose that $\l d(x,\partial\O)\to + \infty$ and $\e \log \l \to 0$ as $\e\to 0$. Then we have the following expansion
\begin{align*}
J_\e(P\d_{x,\l})&= \frac{S_n^{\frac{p-1-\e}{p+1-\e}}}{K(x)^{\frac{2}{p+1-\e}}} \left[ 1- \frac{(n-4)c_2 \D K(x)}{2n^2 S_n K(x) \l^2}\right.\\
&+ \frac{n-4}{n}\e \left(\log \l^{\frac{n-4}{2}} + \frac{c_3}{S_n}\right) + \frac{c_1 H(x,x)}{S_n\l^{n-4}}\\
&+O\left(\frac{\e\log\l}{\l^2} + \frac{1}{\l^{n-3}} + \frac{1}{(\l d)^{n-2}} + \sum_{j=3}^{n-4}\frac{|D^jK(x)|}{\l^j}\right)\\
&\left.+ O\left(\frac{\e \log \l}{(\l d)^{n-4}} + \e^2\log^2 \l +
(\mbox{ if } n < 8) \frac{1}{(\l d)^{2(n-4)}}\right)\right],
\end{align*}
where $S_n$, $c_1$, $c_2$ and $c_3$ are defined in Lemma \ref{l:22}.
\end{pro}
\begin{pf}
According to \cite{BH}, we have
\begin{eqnarray}\label{e:212}
||P\d_{x,\l}||^2 = S_n -c_1 \frac{H(x,x)}{\l^{n-4}} +
O\left(\frac{1}{(\l d)^{n-2}}\right).
\end{eqnarray}
We also have
\begin{align}\label{e:213}
\int_\O K(y)P\d_{x,\l}^{p+1-\e}& = \int_\O K(y)\left(\d_{x,\l} - \var_{x,\l}\right)^{p+1-\e}\notag\\
&= \int_\O K(y)\d_{x,\l}^{p+1-\e} - (p+1-\e)\int_\O K(y)\d_{x,\l}^{p-\e} \var_{x,\l}\notag\\
& + O\left(\int_{ B(x,d)} \d_{x,\l}^{p-1-\e} \var^2_{x,\l}+\frac{1}{(\l d)^{n-1}}\right).
\end{align}
We now observe that, for $n\geq 8$, we have $\frac{n}{n-4}\leq 2$ and thus in this case we have
\begin{align}\label{e:214}
\int_{ B(x,d)} \d_{x,\l}^{p-1-\e} \var^2_{x,\l}& \leq \int_{ B(x,d)} \d_{x,\l}^{\frac{n}{n-4}} \var_{x,\l}^{\frac{n}{n-4}-\e}\notag\\
&\leq ||\var_{x,\l}||_\infty^{\frac{n}{n-4}-\e}\int_{B(x,d)} \d_{x,\l}^{\frac{n}{n-4}}\notag\\
&= O\left(\frac{(\l d^2)^{\e \frac{n-4}{2}} \log (\l d)}{(\l
d)^n}\right)= O\left(\frac{1}{(\l d)^{n-1}}\right),
\end{align}
and, for $n< 8$, we have
\begin{eqnarray}\label{e:215}
\int_{ B(x,d)} \d_{x,\l}^{p-1-\e} \var^2_{x,\l} \leq \frac{1}{(\l d^2)^{(n-4)}} \int_{B(x,d)} \d_{x,\l}^{p-1-\e} = O\left( \frac{1}{(\l d)^{2(n-4)}}\right).
\end{eqnarray}
Thus using Lemma \ref{l:22} and \eqref{e:214} , \eqref{e:215}, we obtain
\begin{align}\label{e:216}
\int_\O K(y)P\d_{x,\l}^{p+1-\e} & = K(x)S_n+\frac{c_2\D K(x)}{2n\l ^2}-\e K(x)S_n\left( \log \l ^{\frac{n-4}{2}}+\frac{c_3}{S_n}\right)-\frac{c_1 2n
K(x)H(x,x)}{(n-4)\l^{n-4}}\notag\\
&+ O\left( \frac{\e \log \l}{(\l d)^{n-4}} + \frac{1}{(\l d)^{n-2}}+ \sum_{j=3}^{n-4} \frac{|D^jK(x)|}{\l ^j}+\frac{1}{\l ^{n-3}}\right)\notag\\
& +O\left(\frac{\e \log \l}{\l ^2}+(\e \log \l )^2 +(if n<8)\frac{1}{(\l d)^{2(n-4)}}\right).
\end{align}

\eqref{e:212}, \eqref{e:216} obviously show that Proposition \ref{p:24} holds.
\end{pf}

The following lemma gives the  basic property of the functional $l_\e$ defined in \eqref{e:23}.
\begin{lem}\label{l:25}
Assume that $x\in\O$ such that $d=d(x,\partial\O) \geq d_0 >0$, and let $v_\e$ be the function obtained
in Proposition \ref{p:21}. Then the functional $l_\e$ has the following expansion :
$$
l_\e (P\d_{x,\l} + v_\e) = \frac{1}{K(x)}\left[ 1 + O\left(
\frac{1}{\l^{n-4}} + \e \log \l + \sum_{j=2}^{n-4}
\frac{|D^jK(x)|}{\l^j}+\sum_{j=1}^k\frac{|D^jK(x)|^2}{\l ^{2j}}+\frac{1}{\l ^{2k+2}}
\right)\right],
$$
where $k$ is the biggest positive integer satisfying $k\leq \frac{n-4}{2}$ .
\end{lem}
\begin{pf}
We have
\begin{eqnarray}\label{e:217}
||P\d_{x,\l} + v_\e||^2 = ||P\d_{x,\l}||^2 + ||v_\e||^2.
\end{eqnarray}
We also have
\begin{eqnarray}\label{e:218}
\int_\O K(y) |P\d_{x,\l} + v_\e|^{p+1-\e} = \int_\O K(y) P\d_{x,\l}^{p+1-\e} + (p+1-\e) \int_\O K(y) P\d_{x,\l}^{p-\e} v_\e + O(||v_\e||^2).
\end{eqnarray}
Thus, using \eqref{e:212}, \eqref{e:216}, \eqref{e:217}, \eqref{e:218} Lemma \ref{l:23} and Proposition \ref{p:21}, we easily derive our lemma.
\end{pf}
\begin{lem}\label{l:26}
Assume that $x\in\O$ such that $d=d(x,\partial\O) \geq d_0 >0$, and let $v_\e$ be the function obtained in Proposition \ref{p:21}. Then
 the following expansion holds.
\begin{align*}
\biggl(\n J_\e (P\d_{x,\l} + v_\e)&, \frac{\partial P\d_{x,\l}}{\partial\l} \biggr)
= \frac{1}{\left(S_n K(x)\right)^{\frac{2}{p+1-\e}}}\left[
\frac{c_2(n-4) \D K(x)}{n^2 K(x) \l^3}-\frac{c_1(n-4)H(x,x)}{\l^{n-3}}\right.\\
& +\frac{(n-4)^2S_n \e}{2n \l} + O\left( \frac{\e \log \l}{\l^3} + \frac{1}{\l^{n-2}}+
\frac{\e ^2\log \l}{\l}+\frac{\e \log \l}{\l^{n-3}} \right)\\
& \left. + O\left( \sum_{3}^{n-4}\frac{|D^jK(x)|}{\l^{j+1}} +
\sum_1^k \frac{|D^jK(x)|^2}{\l^{2j+1}}+\frac{1}{\l
^{2k+3}}+\left(if n<8 \right) \frac{1}{\l ^{2n-7}}  \right) \right],
\end{align*}
where $k$ is the biggest positive integer satisfying $k\leq \frac{n-4}{2}$.
\end{lem}
\begin{pf}
We have
\begin{align}\label{e:219}
\left(\n J_\e (P\d_{x,\l} + v_\e), \frac{\partial P\d_{x,\l}}{\partial\l} \right)&= \frac{2}{\left(\int_\O K(y) |P\d_{x,\l} + v_\e|^{p+1-\e}\right)^{\frac{2}{p+1-\e}}} \left[ \left(P\d_{x,\l}, \frac{\partial P\d_{x,\l}}{\partial\l} \right)\right.\notag\\
&\left. - l_\e(P\d_{x,\l} + v_\e) \int_\O K(y) |P\d_{x,\l}+v_\e|^{p-\e}\frac{\partial P\d_{x,\l}}{\partial\l}\right].
\end{align}
According to \cite{BH}, we have
\begin{eqnarray}\label{e:220}
\left(P\d_{x,\l}, \frac{\partial P\d_{x,\l}}{\partial\l} \right) = \frac{c_1(n-4)H(x,x)}{2 \l^{n-3}} + O \left(\frac{1}{\l^{n-1}}\right).
\end{eqnarray}
On the other hand it follows from Lemma \ref{l:23} and Proposition \ref{p:21} that
\begin{align}\label{e:221}
\int_\O K(y)& |P\d_{x,\l}+ v_\e|^{p-\e}\frac{\partial P\d_{x,\l}}{\partial\l}= \int_\O K(y) P\d_{x,\l}^{p-\e}\frac{\partial P\d_{x,\l}}{\partial\l}\notag\\
& + (p-\e) \int_\O K(y) P\d_{x,\l}^{p-1-\e}v_\e\frac{\partial P\d_{x,\l}}{\partial\l} + O\left(\frac{||v_\e||^2}{\l}\right)\notag\\
&=\int_\O K(y) P\d_{x,\l}^{p-\e}\frac{\partial
P\d_{x,\l}}{\partial\l} + O\left(\sum_1^k \frac{|D^jK(x)|^2}{\l^{2j+1}} +\frac{1}{\l^{2k+3}}+ \frac{1}{\l^{n-3+2\th}} + \frac{\e ^2}{\l}\right).
\end{align}
We are now going to estimate the integral in the right-hand side of \eqref{e:221}.  To this aim, we write
\begin{align}\label{e:222}
\int_\O K(y) P\d_{x,\l}^{p-\e}\frac{\partial P\d_{x,\l}}{\partial\l}&=\int_\O K(y) (\d_{x,\l}- \var_{x,\l})^{p-\e}\frac{\partial (\d_{x,\l}-\var_{x,\l})}{\partial\l}\notag\\
&=\int_\O K(y) \d_{x,\l}^{p-\e}\frac{\partial \d_{x,\l}}{\partial\l} -
\int_\O K(y) \d_{x,\l}^{p-\e}\frac{\partial \var_{x,\l}}{\partial\l}\notag\\
&-(p-\e)\int_\O K(y) \d_{x,\l}^{p-1-\e}\var_{x,\l}\frac{\partial \d_{x,\l}}{\partial\l}
+ O\left(\int_\O \d_{x,\l}^{p-1-\e}\var_{x,\l}|\frac{\partial \var_{x,\l}}{\partial\l}|\right)\notag\\
& +  O\left(\int_\O \d_{x,\l}^{p-1-\e}\frac{\var_{x,\l}^2}{\l}+\frac{1}{\l^{n-1}}\right).
\end{align}
As in \eqref{e:214}  and\eqref{e:215} we derive that
\begin{align}\label{e:223}
O\left(\int_\O \d_{x,\l}^{p-1-\e}\var_{x,\l}|\frac{\partial
\var_{x,\l}}{\partial\l}|\right) =  O\left(\int_\O
\d_{x,\l}^{p-1-\e}\frac{\var_{x,\l}^2}{\l}\right)=O\left(
\frac{1}{\l ^n}+\left(if n<8\right)\frac{1}{\l ^{2n-7}}\right).
\end{align}
Lemma \ref{l:26} follows from \eqref{e:219},..., \eqref{e:223} and
Lemmas  \ref{l:25}, \ref{l:22}.
\end{pf}
\begin{lem}\label{l:27}
Suppose that $K$ satisfies the assumptions of Theorem \ref{t:14} and
\begin{eqnarray}\label{e:229}
|x-x_0| \leq \e^{1/L},\quad \l\in [ c\e^{-1/(n-4)}, c' \e^{-1/(n-4)}].
\end{eqnarray}
Then
$$
||v_\e|| = O\left(\e^{(1+\sigma)/2}\right),
$$
where $\sigma$ is a positive constant and where $v_\e$ is defined in Proposition \ref{p:21}.
\end{lem}
\begin{pf}
It view of Proposition \ref{p:21}, we only need to check
\begin{eqnarray}\label{e:230}
\frac{|D^jK(x)|}{\l^j} = O(\e^{1+\sigma}).
\end{eqnarray}
But, by assumptions imposed on $K$, we see that if $\sigma >0$ is small enough, then
$$
\frac{|D^jK(x)|}{\l^j} \leq C \frac{|x-x_0|^{L-j}}{\l^j}\leq C \e^{\frac{L-j}{L}} \e^{\frac{j}{n-4}}= O\left(\e^{1+\sigma}\right).
$$
\end{pf}
\begin{lem}\label{l:28}
Suppose that $K$ satisfies the assumptions of Theorem \ref{t:14} and \eqref{e:229} holds. Then we have the following estimates:
\begin{eqnarray*}
\begin{array}{cccc}
1.\qquad\qquad &\left(\n J_\e(P\d_{x,\l} + v_\e), P\d_{x,\l}\right)&=&O\left(\e^{1-\sigma}\right)\quad\qquad\\
2.\qquad\qquad &\left(\n J_\e(P\d_{x,\l} + v_\e), \frac{\partial P\d_{x,\l}}{\partial\l}\right)&=&O\left(\e^{1+ \frac{1}{n-4}}\right)\,\,\quad\,\\
3.\qquad\qquad &\left(\n J_\e(P\d_{x,\l} + v_\e), \frac{\partial P\d_{x,\l}}{\partial x_j}\right)&=&O\left(\e^{1+\sigma - \frac{1}{n-4}}\right),
\end{array}
\end{eqnarray*}
where $v_\e$ is defined in Proposition \ref{p:21}.
\end{lem}
\begin{pf}
Lemma \ref{l:26} and \eqref{e:230} give Claim 2. To prove Claim 1, we write
\begin{align*}
\left(\n J_\e (P\d_{x,\l} + v_\e), P\d_{x,\l} \right)&= \frac{2}{\left(\int_\O K(y) |P\d_{x,\l} + v_\e|^{p+1-\e}\right)^{\frac{2}{p+1-\e}}} \left[ \left(P\d_{x,\l},  P\d_{x,\l} \right)\right.\\
&\left. - l_\e(P\d_{x,\l} + v_\e) \int_\O K(y) |P\d_{x,\l}+v_\e|^{p-\e} P\d_{x,\l}\right]
\end{align*}
and thus, using Lemmas \ref{l:23}, \ref{l:25},  Proposition \ref{p:21}, \eqref{e:212}, \eqref{e:216} and \eqref{e:230}\ we easily derive Claim 1.\\
As in \eqref{e:220}, \eqref{e:222} (see also \cite{BH})we have
\begin{eqnarray}\label{e:231}
\left(P\d_{x,\l}, \frac{\partial P\d_{x,\l}}{\partial x_i} \right) = \frac{\partial H(x,x)}{\partial x_i} \frac{c_1}{2 \l ^{n-4}} + O \left(\frac{1}{\l ^{n-2}}\right).
\end{eqnarray}
 \begin{align}\label{e:232}
\int_\O K(y) P\d_{x,\l}^{p-\e}\frac{\partial P\d_{x,\l}}{\partial
x_i}= O\left(\sum_1^{n-4} \frac{|D^jK(x)|}{\l^{j-1}}
+\frac{1}{\l^{n-4}}+(if n<8)\frac{1}{\l ^{2n-9}}\right).
\end{align}
 Then Claim 3 follows.
\end{pf}

 Next, our goal is to estimate $||\partial v_\e/\partial\l||$, where $v_\e$ is defined in Proposition \ref{p:21}. To this aim, we follow \cite{CY}, namely, we write the following
 decomposition
 \begin{eqnarray} \label{e:233}
 \frac{\partial v_\e}{\partial\l}= w + \a P\d_{x,\l} + \b \frac{\partial P\d_{x,\l}}{\partial\l} + \sum_{j=1}^n \g_j \frac{\partial P\d_{x,\l}}{\partial x_j},
\end{eqnarray}
 where $\a$, $\b$ and $\g_j$ are chosen in such a way that $w \in E_{x,\l}$.
\begin{lem}\label{l:29}
Let $\a$, $\b$ and $\g_j$ be coefficients in \eqref{e:233} and assume that \eqref{e:229} holds. Then we have the following estimates
$$
\a = O\left(\frac{||v_\e||}{\l}\right),\quad \b =O(||v_\e||), \quad \g_j = O\left(\frac{||v_\e||}{\l^2}\right).
$$
\end{lem}
\begin{pf}
Taking the scalar product of \eqref{e:233} with $P\d_{x,\l}$, $\partial P\d_{x,\l}/\partial\l$ and $\partial P\d_{x,\l}/\partial x_i$ for $i=1,...,n$, we obtain
\begin{eqnarray*}
\a ||P\d_{x,\l}||^2 + \b \left(\frac{\partial P\d_{x,\l}}{\partial\l}, P\d_{x,\l}\right) + \sum_{j=1}^n \g_j \left(\frac{\partial P\d_{x,\l}}{\partial x_j}, P\d_{x,\l}\right) = 0,
\end{eqnarray*}
\begin{align*}
\a \left(P\d_{x,\l}, \frac{\partial P\d_{x,\l}}{\partial\l}\right)& + \b \biggl|\biggl|\frac{\partial P\d_{x,\l}}{\partial\l}\biggr|\biggr|^2 +
\sum_{j=1}^n \g_j \left(\frac{\partial P\d_{x,\l}}{\partial x_j},\frac{\partial P\d_{x,\l}}{\partial\l}\right)\\
&= -\left(v_\e, \frac{\partial^2 P\d_{x,\l}}{\partial\l^2}\right)= O\left(\frac{||v_\e||}{\l^2}\right),
\end{align*}
\begin{align*}
\a \left(P\d_{x,\l}, \frac{\partial P\d_{x,\l}}{\partial x_i}\right)& +
\b \left(\frac{\partial P\d_{x,\l}}{\partial\l}, \frac{\partial P\d_{x,\l}}{\partial x_i}\right) +
\sum_{j=1}^n \g_j \left(\frac{\partial P\d_{x,\l}}{\partial x_j},\frac{\partial P\d_{x,\l}}{\partial x_i}\right)\\
&=-\left(v_\e, \frac{\partial^2 P\d_{x,\l}}{\partial\l\partial x_i}\right)= O\left(||v_\e||\right).
\end{align*}
Thus, we derive that
\begin{align*}
\frac{\a}{\l}(S_n + O(\e)) + \frac{\b}{\l^2}O(\e) + \sum_1^n \g_j O(\e)&=0\\
\frac{\a}{\l}O(\e)+ \frac{\b}{\l^2}(c_n' + O(\e)) + \sum_1^n \g_j O(\e)&= O\left(\frac{||v_\e||}{\l^2}\right)\\
\frac{\a}{\l}O(\e) + \frac{\b}{\l^2}O(\e) + \sum_{j\ne i}\g_j O(\e)+ \g_i(c_n^{''}+O(\e))&= O\left(\frac{||v_\e||}{\l^2}\right).
\end{align*}
Solving the above system we get the desired estimates.
\end{pf}

Now, for a fixed $w_0\in E_{x_0,\l_0}$, we denote $\pi(x,\l)$ the orthogonal projection of $w_0$ onto $E_{x,\l}$.
We then have
\begin{eqnarray}\label{e:234}
w_0= \pi(x,\l) + a(x,\l)P\d_{x,\l} + b(x,\l)\frac{\partial P\d_{x,\l}}{\partial\l}+ \sum_{j=1}^n g_j(x,\l)\frac{\partial P\d_{x,\l}}{\partial x_j}.
\end{eqnarray}
\begin{lem}\label{l:210}
The map $\pi(.,.)$ is $C^1$ with respect to $x$ and $\l$, and
\begin{align*}
a(x_0,\l_0)&=0,\qquad \frac{\partial a(x_0,\l_0)}{\partial\l}= O\left(\frac{||w_0||}{\l}\right),\\
b(x_0,\l_0)&=0,\qquad \frac{\partial b(x_0,\l_0)}{\partial\l}= O\left(||w_0||\right),\\
g_j(x_0,\l_0)&=0,\qquad \frac{\partial g_j(x_0,\l_0)}{\partial\l}= O\left(\frac{||w_0||}{\l^2}\right).
\end{align*}
\end{lem}
\begin{pf}
First of all, we easily deduce from the fact that $w_0 \in E_{x_0,\l_0}$ the following:
$$
a(x_0,\l_0)=b(x_0,\l_0)=g_j(x_0,\l_0)=0.
$$
Secondly, it is clear to see that $a(x,\l)$, $b(x,\l)$ and $g_j(x,\l)$ satisfy
\begin{eqnarray}\label{e:235}
a||P\d_{x,\l}||^2 + b \left(\frac{\partial P\d_{x,\l}}{\partial\l}, P\d_{x,\l}\right) + \sum_{j=1}^n g_j \left(\frac{\partial P\d_{x,\l}}{\partial x_j}, P\d_{x,\l}\right) = \left(w_0, P\d_{x,\l}\right),
\end{eqnarray}
\begin{eqnarray}\label{e:236}
a \left(P\d_{x,\l}, \frac{\partial P\d_{x,\l}}{\partial\l}\right) +
b \biggl|\biggl|\frac{\partial P\d_{x,\l}}{\partial\l}\biggr|\biggr|^2 +
\sum_{j=1}^n g_j \left(\frac{\partial P\d_{x,\l}}{\partial x_j},\frac{\partial P\d_{x,\l}}{\partial\l}\right)= \left(w_0,\frac{\partial P\d_{x,\l}}{\partial\l}\right),
\end{eqnarray}
\begin{eqnarray}\label{e:237}
a \left(P\d_{x,\l}, \frac{\partial P\d_{x,\l}}{\partial x_i}\right) + b \left(\frac{\partial P\d_{x,\l}}{\partial\l}, \frac{\partial P\d_{x,\l}}{\partial x_i}\right) + \sum_{j=1}^n g_j \left(\frac{\partial P\d_{x,\l}}{\partial x_j},\frac{\partial P\d_{x,\l}}{\partial x_i}\right)= \left(w_0, \frac{\partial P\d_{x,\l}}{\partial x_i}\right).
\end{eqnarray}
Solving the above system we easily see that  $a(x,\l)$, $b(x,\l)$ and $g_j(x,\l)$ are $C^1$ with respect to $x$ and $\l$.
Differentiating \eqref{e:235}, \eqref{e:236} and \eqref{e:237} with respect to $\l$, we obtain
\begin{align*}
\frac{\partial a(x_0,\l_0)}{\partial\l}||P\d_{x_0,\l_0}||^2 &+ \frac{\partial b(x_0,\l_0)}{\partial\l} \left(\frac{\partial P\d_{x_0,\l_0}}{\partial\l}, P\d_{x_0,\l_0}\right) + \sum_{j=1}^n \frac{\partial g_j(x_0,\l_0)}{\partial\l}\\
&\times \left(\frac{\partial P\d_{x_0,\l_0}}{\partial x_j}, P\d_{x_0,\l_0}\right)
 = \left(w_0, \frac{\partial P\d_{x_0,\l_0}}{\partial\l}\right)=0,
\end{align*}
\begin{align*}
 \frac{\partial a(x_0,\l_0)}{\partial\l}&\left(P\d_{x_0,\l_0}, \frac{\partial P\d_{x_0,\l_0}}{\partial\l}\right) +
 \frac{\partial b(x_0,\l_0)}{\partial\l}  \biggl|\biggl|\frac{\partial P\d_{x_0,\l_0}}{\partial\l}\biggr|\biggr|^2 +
  \sum_{j=1}^n \frac{\partial g_j(x_0,\l_0)}{\partial\l}\\
&\times \left(\frac{\partial P\d_{x_0,\l_0}}{\partial x_j},\frac{\partial P\d_{x_0,\l_0}}{\partial\l}\right)
= \left(w_0,\frac{\partial^2 P\d_{x_0,\l_0}}{\partial\l^2}\right)= O\left(\frac{||w_0||}{\l^2}\right),
\end{align*}
\begin{align*}
\frac{\partial a(x_0,\l_0)}{\partial\l}&\left(P\d_{x_0,\l_0}, \frac{\partial P\d_{x_0,\l_0}}{\partial x_i}\right) +
 \frac{\partial b(x_0,\l_0)}{\partial\l} \left(\frac{\partial P\d_{x_0,\l_0}}{\partial\l}, \frac{\partial P\d_{x_0,\l_0}}{\partial x_i}\right) \\
&+ \sum_{j=1}^n \frac{\partial g_j(x_0,\l_0)}{\partial\l} \left(\frac{\partial P\d_{x_0,\l_0}}{\partial x_j},\frac{\partial P\d_{x_0,\l_0}}{\partial x_i}\right)
= \left(w_0, \frac{\partial^2 P\d_{x_0,\l_0}}{\partial\l \partial x_i}\right)= O(||w_0||).
\end{align*}
Thus, as in the proof of Lemma \ref{l:29}, we derive the desired result.
\end{pf}
\begin{pro}\label{p:211}
Assume that \eqref{e:229} holds. Then, we have the following estimate
$$
\biggl|\biggl|\frac{\partial v_\e}{\partial\l}\biggr|\biggr| = O \left(\e^{\frac{1+\sigma}{2} + \frac{1}{n-4}}\right),
$$
where $v_\e$ is defined in Proposition \ref{p:21}.
\end{pro}
\begin{pf}
In view of \eqref{e:233} and Lemma \ref{l:29}, we only need to estimate $||w||$. Let $\pi(x',\l')$ be the orthogonal projection of $w\in E_{x,\l}$ onto $E_{x',\l'}$. Thus we have
\begin{eqnarray}\label{e:238}
\left(\n J_\e \left( P\d_{x',\l'} + v_\e(x',\l')\right), \pi(x',\l')\right)= 0.
\end{eqnarray}
Differentiating \eqref{e:238} with respect to $\l'$ and letting $(x',\l')=(x,\l)$, we obtain
\begin{eqnarray}\label{e:239}
D^2J_\e(P\d_{x,\l}+ v_\e)\left(\frac{\partial P\d_{x,\l}}{\partial\l}+ \frac{\partial v_\e}{\partial\l}, w\right) + \left(\n J_\e (P\d_{x,\l}+ v_\e), \frac{\partial \pi (x,\l)}{\partial\l}\right)=0.
\end{eqnarray}
It follows from Lemmas \ref{l:28} and \ref{l:210} that
\begin{align}\label{e:240}
\biggl(\n J_\e & (P\d_{x,\l}+ v_\e), \frac{\partial \pi (x,\l)}{\partial\l}\biggr)= \frac{\partial a}{\partial\l} \left(\n J_\e(P\d_{x,\l}+ v_\e), P\d_{x,\l}\right)\notag\\
&+ \frac{\partial b}{\partial\l} \left(\n J_\e(P\d_{x,\l}+ v_\e),\frac{\partial P\d_{x,\l}}{\partial\l}\right)+  \sum_1^n \frac{\partial g_j}{\partial\l} \left(\n J_\e(P\d_{x,\l}+ v_\e),\frac{\partial P\d_{x,\l}}{\partial x_j}\right)\notag\\
&=O\left(||w||\left(\frac{\e^{1-\sigma}}{\l} + \e^{1-\sigma + \frac{1}{n-4}} + \frac{ \e^{1-\sigma - \frac{1}{n-4}}}{\l^2}\right)\right)= O\left(||w|| \e^{1-\sigma + \frac{1}{n-4}}\right).
\end{align}
Combining \eqref{e:238} and \eqref{e:239} and taking Lemmas \ref{l:28} and \ref{l:29} into account we obtain
\begin{align}\label{e:241}
D^2J_\e(&P\d_{x,\l}+ v_\e)(w,w)= -D^2J_\e(P\d_{x,\l}+ v_\e)\biggl(\frac{\partial P\d_{x,\l}}{\partial\l} + \a P\d_{x,\l}+ \b\frac{\partial P\d_{x,\l}}{\partial\l}\notag\\
& + \sum_1^n \g_j \frac{\partial P\d_{x,\l}}{\partial x_j}, w\biggr) + O\left(||w|| \e^{1-\sigma +\frac{1}{n-4}}\right)\notag\\
& = -D^2J_\e(P\d_{x,\l}+ v_\e)\left(\frac{\partial P\d_{x,\l}}{\partial\l}, w\right) + O\left(\frac{||w||||v_\e||}{\l}\right)+  O\left(||w|| \e^{1-\sigma +\frac{1}{n-4}}\right)\notag\\
& = -D^2J_\e(P\d_{x,\l}+ v_\e)\left(\frac{\partial P\d_{x,\l}}{\partial\l}, w\right) + O\left(||w|| \e^{\frac{1+\sigma}{2} +\frac{1}{n-4}}\right).
\end{align}
We now claim that
\begin{eqnarray}\label{e:242}
D^2J_\e(P\d_{x,\l}+ v_\e)(w,w) \geq \rho ||w||^2,
\end{eqnarray}
for some positive constant $\rho$ and
\begin{eqnarray}\label{e:243}
D^2J_\e(P\d_{x,\l}+ v_\e)\left(\frac{\partial P\d_{x,\l}}{\partial\l}, w\right) = O\left(||w|| \e^{\frac{1+\sigma}{2}+ \frac{1}{n-4}}\right).
\end{eqnarray}
Then obviously \eqref{e:241}, \eqref{e:242} and \eqref{e:243} imply that $||w||= O\left(\e^{\frac{1+\sigma}{2}+ \frac{1}{n-4}}\right)$, and Proposition \ref{p:211}  follows.

It remains to prove \eqref{e:242} and \eqref{e:243}. To this aim, we write
\begin{align}\label{e:244}
D^2J_\e&(P\d_{x,\l}+ v_\e)(\var, \psi)= \frac{2(\var, \psi)}{\left(\int_\O K(y)|u|^{p+1-\e}\right)^{\frac{2}{p+1-\e}}}-\frac{4(P\d_{x,\l}+ v_\e, \var)}{\left(\int_\O K(y)|u|^{p+1-\e}\right)^{\frac{2}{p+1-\e}+1}}\notag\\
&\times \int_\O K(y)|P\d_{x,\l}+ v_\e|^{p-\e}\psi
-\frac{4(P\d_{x,\l}+ v_\e, \psi)}{\left(\int_\O K(y)|u|^{p+1-\e}\right)^{\frac{2}{p+1-\e}+1}}\int_\O K(y)|P\d_{x,\l}+ v_\e|^{p-\e}\var\notag\\
& +2(p+3-\e)||P\d_{x,\l}+ v_\e||^2 \frac{\int_\O K(y)|P\d_{x,\l}+ v_\e|^{p-\e}\var\int_\O K(y)|P\d_{x,\l}+ v_\e|^{p-\e}\psi}{\left(\int_\O K(y)|P\d_{x,\l}+ v_\e|^{p+1-\e}\right)^{\frac{2}{p+1-\e}+2}}\notag\\
& -\frac{2(p-\e)||P\d_{x,\l}+ v_\e||^2}{\left(\int_\O K(y)|P\d_{x,\l}+ v_\e|^{p+1-\e}\right)^{\frac{2}{p+1-\e}+1}} \int_\O K(y)|P\d_{x,\l}+ v_\e|^{p-1-\e}\var\psi.
\end{align}

{\it Verification of} \eqref{e:242}. First, we notice that
\begin{eqnarray}\label{e:245}
\left(P\d_{x,\l}+ v_\e, w\right)=\left(v_\e, w\right)= O\left( \e^{(1+\sigma)/2}||w||\right),
\end{eqnarray}
where we have used Lemma \ref{l:27}. By Lemma \ref{l:23}, we see
\begin{eqnarray}\label{e:246}
\int_\O K(y)|P\d_{x,\l}+ v_\e|^{p-\e}w= \int_\O K(y)|P\d_{x,\l}|^{p-\e}w + O\left(||v_\e||||w||\right)= O\left(\e^{\frac{1+\sigma}{2}}||w||\right).
\end{eqnarray}
Thus \eqref{e:242} follows from \eqref{e:244},...,\eqref{e:246} and Proposition 3.4 in \cite{BE1}.

{\it Verification of} \eqref{e:243}. We have
\begin{align}
&\left(\frac{\partial P\d_{x,\l}}{\partial\l}, w\right)=0,\label{e:247}\\
&\left(P\d_{x,\l}+ v_\e,\frac{\partial P\d_{x,\l}}{\partial\l}\right)=
 \left(P\d_{x,\l},\frac{\partial P\d_{x,\l}}{\partial\l}\right) = O\left(\frac{1}{\l^{n-3}}\right)=
 O\left(\e^{1+ \frac{1}{n-4}}\right).\label{e:248}
\end{align}
Also, it follows from\eqref{e:221} that
\begin{eqnarray}\label{e:249}
\int_\O K(y)|P\d_{x,\l}+ v_\e|^{p-\e}\frac{\partial P\d_{x,\l}}{\partial\l}= O\left(\e^{1+ \frac{1}{n-4}}\right),
\end{eqnarray}
and  by Lemma \ref{l:23} as in \eqref{e:246}, we have
\begin{align} \label{e:250}
&\ \int_\O K(y)|P\d_{x,\l}+ v_\e|^{p-1-\e}\frac{\partial P\d_{x,\l}}{\partial\l}w=
\int_\O K(y)|P\d_{x,\l}|^{p-1-\e}\frac{\partial P\d_{x,\l}}{\partial\l}w +
O\left(\frac{||v_\e||||w||}{\l}\right)\notag \\
&\ =O\left(||w||\e^{\frac{1+\sigma}{2}+\frac{1}{n-4}}\right).
\end{align}
Combining \eqref{e:245},..., \eqref{e:250} we obtain \eqref{e:243} and this completes the proof of Proposition \ref{p:211}.
\end{pf}
\begin{lem}\label{l:212}
The derivative of the functional $J_\e$ satisfies
\begin{align*}
(i)\quad \frac{\partial}{\partial\l} \biggl (\n J_\e(P\d_{x,\l}+ v_\e), \frac{\partial P\d_{x,\l}}{\partial\l}\biggr )&=\frac{1}{\left(K(x)S_n\right)^{\frac{2}{p+1-\e}}}\biggl [ -\frac{(n-4)^2S_n\e}{2n\l^2}\\
&+ \frac{c_1(n-4)(n-3)H(x,x)}{\l^{n-2}} + O\left(\e^{1+\sigma + 2/(n-4)}\right)\biggr ],\\
(ii)\quad \frac{\partial}{\partial\l} \biggl (\n J_\e(P\d_{x,\l}+ v_\e), \frac{\partial P\d_{x,\l}}{\partial x_j}\biggr )&= O\left(\e^{1-\sigma}\right),\\
(iii)\qquad \frac{\partial}{\partial\l} \left(\n J_\e(P\d_{x,\l}+ v_\e), P\d_{x,\l}\right)&= O\left(\e^{1 + 1/(n-4)}\right).
\end{align*}
\end{lem}
\begin{pf}
By easy computations we have
\begin{align*}
\frac{\partial}{\partial\l} \biggl (\n J_\e(P\d_{x,\l}+ v_\e), \frac{\partial P\d_{x,\l}}{\partial\l}\biggr )&= D^2J_\e(P\d_{x,\l} + v_\e)\left(\frac{\partial P\d_{x,\l}}{\partial\l} + \frac{\partial v_\e}{\partial\l}, \frac{\partial P\d_{x,\l}}{\partial\l}\right)\\
&+ \left(\n J_\e(P\d_{x,\l} + v_\e), \frac{\partial^2 P\d_{x,\l}}{\partial\l^2}\right).
\end{align*}
First, we estimate $D^2J_\e(P\d_{x,\l} + v_\e)\left(\frac{\partial P\d_{x,\l}}{\partial\l}, \frac{\partial v_\e}{\partial\l}\right)$. Using Proposition \ref{p:211}, we obtain
\begin{align} \label{e:251}
\left(P\d_{x,  \l} + v_\e, \frac{\partial v_\e}{\partial\l}\right)&= \left(v_\e, \frac{\partial v_\e}{\partial\l}\right)= O\left(\e^{1+\sigma + \frac{1}{n-4}}\right).
\end{align}
As in the proof of Lemma \ref{l:26}, we have
\begin{align}\label{e:252}
\int_\O K(y)|P\d_{x,\l} + v_\e|^{p-\e}\frac{\partial v_\e}{\partial\l}&= O\left(\e^{(1 + \sigma)/2 + 1/(n-4)}\right).
\end{align}
Combining \eqref{e:248}, \eqref{e:249}, \eqref{e:251}, \eqref{e:252}, we obtain
\begin{align}\label{e:253}
D^2&J_\e(P\d_{x,\l} + v_\e)\left(\frac{\partial P\d_{x,\l}}{\partial\l}, \frac{\partial v_\e}{\partial\l}\right)= \frac{2}{\left(\int_\O K(y)|P\d_{x,\l}+ v_\e|^{p+1-\e}\right)^{\frac{2}{p+1-\e}}} \biggl [ \left(\frac{\partial P\d_{x,\l}}{\partial\l}, \frac{\partial v_\e}{\partial\l}\right)\notag\\
&- (p-\e)l_\e(P\d_{x,\l}+v_\e)\int_\O K(y)|P\d_{x,\l} + v_\e|^{p-1-\e} \frac{\partial P\d_{x,\l}}{\partial\l} \frac{\partial v_\e}{\partial\l}\biggr ] + O\left(\e^{1+\sigma + \frac{2}{n-4}}\right).
\end{align}
We now notice that
\begin{align}
&\frac{\partial v_\e}{\partial\l}= w +\a P\d_{x,\l} +\b \frac{\partial P\d_{x,\l}}{\partial\l} + \sum_{j=1}^n \g_j \frac{\partial P\d_{x,\l}}{\partial x_j},\label{e:254}\\
&\biggl(\frac{\partial P\d_{x,\l}}{\partial\l}, w\biggr)= 0.\label{e:255}
\end{align}
Consequently, it follows from Lemma \ref{l:23} and Proposition \ref{p:211} that
\be\label{e:256}
\int_\O K(y)|P\d_{x,\l} + v_\e|^{p-1-\e} \frac{\partial P\d_{x,\l}}{\partial\l} w= O\left(\e^{1+\sigma + 2/(n-4)}\right).
\ee
Now, using Lemmas \ref{l:23}, and \ref{l:29}, we obtain
\begin{align}
\a  \int_\O K(y)|P\d_{x,\l} +& v_\e|^{p-1-\e} P\d_{x,\l} \frac{\partial P\d_{x,\l}}{\partial\l} =  O\left(\e^{1+\sigma + 2/(n-4)}\right)\label{e:257}\\
& \a \left(\frac{\partial P\d_{x,\l}}{\partial\l}, P\d_{x,\l}\right)= O\left(\e^{1+\sigma + 2/(n-4)}\right).\label{e:258}
\end{align}
In the same way, we have
\begin{align}
 \g_j  \int_\O K(y)|P\d_{x,\l} + v_\e|^{p-1-\e} \frac{\partial P\d_{x,\l}}{\partial x_j} \frac{\partial
 P\d_{x,\l}}{\partial\l}= O\left(\e^{1+\sigma + 2/(n-4)}\right) \label{e:259}\\
\g_j \left(\frac{\partial P\d_{x,\l}}{\partial\l},\frac{\partial P\d_{x,\l}}{\partial x_i}\right)= O\left(\e^{1+\sigma + 2/(n-4)}\right).\label{e:260}
\end{align}
As in \eqref{e:221}, \eqref{e:222} and using Lemma \ref{l:25} we
have
\begin{align*}
\b (p-\e ) l_\e(P\d_{x,\l}+v_\e)\int_\O K(y)|P\d_{x,\l}
+v_\e|^{p-1-\e} \biggl| \frac{\partial
P\d_{x,\l}}{\partial\l}\biggr|^2=\b (p-\e )l_\e(P\d_{x,\l}+v_\e)
\\ \times \int_\O K(y)P\d_{x,\l} ^{p-1-\e} \biggl|\frac{\partial
P\d_{x,\l}}{\partial\l}\biggr|^2+O\left(\frac{\b ||v_\e||}{\l
^2}\right)=\b \biggl|\biggl|\frac{\partial
P\d_{x,\l}}{\partial\l}\biggl|\biggl|^2 + O\left(\e^{1+\sigma + 2/(n-4)}\right),
\end{align*}
then, by Lemma \ref{l:29}, we derive that
 \begin{align}\label{e:261}
\b \biggr|\biggr|\frac{\partial P\d_{x,\l}}{\partial\l} \biggr|\biggr|^2-\b (p-\e )
l_\e(P\d_{x,\l}+v_\e)\int_\O K(y)|P\d_{x,\l}
 +v_\e|^{p-1-\e} \biggl| \frac{\partial
 P\d_{x,\l}}{\partial\l}\biggr|^2= O\left(\e^{1+\sigma +
 2/(n-4)}\right).
 \end{align}
Combining \eqref{e:251},..., \eqref{e:261} we obtain
\begin{equation}\label{e:259'}
D^2J_\e(P\d_{x,\l} + v_\e)\left(\frac{\partial P\d_{x,\l}}{\partial\l}, \frac{\partial v_\e}{\partial\l}\right)= O\left(\e^{1+\sigma + 2/(n-4)}\right).
\end{equation}
We now write
\begin{align*}
 D^2&J_\e(P\d_{x,\l} + v_\e)\left(\frac{\partial P\d_{x,\l}}{\partial\l}, \frac{\partial P\d_{x,\l}}{\partial\l}\right)
+ \left(\n J_\e(P\d_{x,\l} + v_\e), \frac{\partial^2 P\d_{x,\l}}{\partial\l^2}\right)\notag\\
&=\frac{2}{(\int_\O K(y)|P\d_{x,\l}+v_\e|^{p+1-\e})^{\frac{2}{p+1-\e}}}
 \biggl\{ \biggl [ \biggl |\biggl | \frac{\partial P\d_{x,\l}}{\partial\l}\biggr |\biggr |^2 - (p-\e)l_\e(P\d_{x,\l} + v_\e)\notag\\
&\times \int_\O K(y)|P\d_{x,\l}+v_\e|^{p-1-\e}\biggl |\frac{\partial P\d_{x,\l}}{\partial\l}\biggr |^2 \biggr ]
+\left(P\d_{x,\l} + v_\e, \frac{\partial^2 P\d_{x,\l}}{\partial\l^2}\right)\notag \\
  & -
l_\e(P\d_{x,\l} + v_\e)\int_\O K(y)|P\d_{x,\l}+v_\e|^{p-\e}\frac{\partial^2 P\d_{x,\l}}{\partial\l^2}\biggr \} +O\left(\e^{1+\sigma +
2/(n-4)}\right).\notag\\
\end{align*}
\begin{align}
\qquad & =\frac{2}{(\int_\O K(y)|P\d_{x,\l}+v_\e|^{p+1-\e})^{\frac{2}{p+1-\e}}}
 \biggl\{ \biggl  |\biggl | \frac{\partial P\d_{x,\l}}{\partial\l}\biggl |\biggl |^2+\left(P\d_{x,\l} , \frac{\partial^2 P\d_{x,\l}}{\partial\l^2}\right)\notag\\
& - l_\e(P\d_{x,\l} + v_\e) \left((p-\e)\int_\O K(y) P\d_{x,\l} ^{p-1-\e}\biggl |\frac{\partial P\d_{x,\l}}{\partial\l}\biggl |^2 +  \int_\O K(y)P\d_{x,\l}^{p-\e}\frac{\partial^2 P\d_{x,\l}}{\partial\l^2}\right)\notag\\
& -(p-\e )l_\e (P\d_{x,\l} + v_\e )\biggl( (p-1-\e)\int_\O K(y) P\d_{x,\l} ^{p-2-\e}  \biggl |\frac{\partial P\d_{x,\l}}{\partial\l}\biggr |^2v_\e\notag\\
&+  \int_\O K(y)P\d_{x,\l}^{p-1-\e} \frac{\partial^2
P\d_{x,\l}}{\partial\l^2}v_\e\biggr)
 +\left(v_\e,\frac{\partial^2 P\d_{x,\l}}{\partial\l^2}\right)\biggl\} +O\left(\e^{1+\sigma + 2/(n-4)}\right).\label{e:262}
\end{align}
We now observe that
\begin{align}
\biggl |\biggl | \frac{\partial P\d_{x,\l}}{\partial\l}\biggr |\biggr |^2= \int_\O &\D^2\left(\frac{\partial P\d_{x,\l}}{\partial\l}\right)\frac{\partial P\d_{x,\l}}{\partial\l}= p\int_\O \d_{x,\l}^{p-1}\frac{\partial \d_{x,\l}}{\partial\l}\frac{\partial P\d_{x,\l}}{\partial\l},\label{e:263}\\
\left(P\d_{x,\l}, \frac{\partial^2 P\d_{x,\l}}{\partial\l^2}\right)&= \int_\O \D^2\left(\frac{\partial^2 P\d_{x,\l}}{\partial\l^2}\right)P\d_{x,\l}=p(p-1)\int_\O \d_{x,\l}^{p-2}\left(\frac{\partial \d_{x,\l}}{\partial\l}\right)^2 P\d_{x,\l}\notag\\
& +p\int_\O \d_{x,\l}^{p-1}\frac{\partial^2 \d_{x,\l}}{\partial\l^2}P\d_{x,\l}.\label{e:264}
\end{align}
Thus, using \eqref{e:263}, \eqref{e:264} and Proposition 2.1 of \cite{BH}, we obtain
\begin{equation}\label{e:265}
\biggl |\biggl | \frac{\partial P\d_{x,\l}}{\partial\l}\biggr |\biggr |^2 + \left(P\d_{x,\l}, \frac{\partial^2 P\d_{x,\l}}{\partial\l^2}\right)=-\frac{c_1(n-3)(n-4)H(x,x)}{2\l^{n-2}} + O\left(\frac{1}{\l^n}\right).
\end{equation}
As in the proof of Lemma \ref{l:22} and using  Proposition 2.1 of \cite{BH} , we find
\begin{align}\label{e:266}
(p-\e)
\int_\O K(y) & P\d_{x,\l}^{p-1-\e}\biggl |\frac{\partial P\d_{x,\l}}{\partial\l}\biggr |^2 +\int_\O K(y)P\d_{x,\l}^{p-\e}\frac{\partial^2
P\d_{x,\l}}{\partial\l^2}\notag\\
&=\frac{(n-4)^2 S_n K(x)\e}{4n\l^2} - \frac{c_1(n-4)(n-3)K(x)H(x,x)}{\l^{n-2}} +O\left(\e^{1+\sigma + \frac{2}{n-4}}\right) .
\end{align}
On the other hand, using Lemma \ref{l:25} it is easy to check
\begin{align}\label{e:267}
 - (p-\e)&(p-1-\e)l_\e(P\d_{x,\l} + v_\e)
\int_\O K(y) P\d_{x,\l}^{p-2-\e}\biggl |\frac{\partial P\d_{x,\l}}{\partial\l}\biggr |^2 v_\e\notag\\&
- (p-\e)l_\e(P\d_{x,\l} + v_\e)
\int_\O K(y)P\d_{x,\l}^{p-1-\e}\frac{\partial^2 P\d_{x,\l}}{\partial\l^2} v_\e + \left( v_\e, \frac{\partial^2 P\d_{x,\l}}{\partial\l^2}\right)\notag\\
 &= -\int_\O \frac{\partial^2}{\partial\l^2}\left(\D^2 P\d_{x,\l}\right)v_\e +\left( v_\e, \frac{\partial^2 P\d_{x,\l}}{\partial\l^2}\right) +O\left(\e^{1+\sigma + \frac{2}{n-4}}\right)\notag\\
&= O\left(\e^{1+\sigma + \frac{2}{n-4}}\right).
\end{align}
Combining the above estimates, Claim $(i)$ follows. The proof of Claim $(ii)$ is similar to that of Claim $(i)$ and therefore is omitted. To prove Claim $(iii)$, we write
\begin{align}\label{e:268}
\frac{\partial}{\partial\l} &\biggl (\n J_\e(P\d_{x,\l}+ v_\e), P\d_{x,\l}\biggr )\notag\\
&= D^2J_\e(P\d_{x,\l} + v_\e)\left(\frac{\partial P\d_{x,\l}}{\partial\l} + \frac{\partial v_\e}{\partial\l}, P\d_{x,\l}\right)
+ \left(\n J_\e(P\d_{x,\l} + v_\e), \frac{\partial P\d_{x,\l}}{\partial\l}\right)\notag\\
&= D^2J_\e(P\d_{x,\l} + v_\e)\left(\frac{\partial P\d_{x,\l}}{\partial\l} + \frac{\partial v_\e}{\partial\l}, P\d_{x,\l}\right) + O\left(\e^{1+ 1/(n-4)}\right),
\end{align}
where we have used Lemma \ref{l:28}.\\
Now, as in the proof of \eqref{e:259'}, we obtain
\begin{equation}\label{e:269}
D^2J_\e(P\d_{x,\l} + v_\e)\left( \frac{\partial v_\e}{\partial\l}, P\d_{x,\l}\right) = O\left(\e^{1+\sigma + 1/(n-4)}\right).
\end{equation}
Computations similar to that in the proof of Claim $(i)$ show that
\begin{equation}\label{e:270}
D^2J_\e(P\d_{x,\l} + v_\e)\left( \frac{\partial P\d_{x,\l}}{\partial\l}, P\d_{x,\l}\right) = O\left(\e^{1 + 1/(n-4)}\right).
\end{equation}
Hence, Claim $(iii)$ follows from \eqref{e:268}, \eqref{e:269} and \eqref{e:270}. This completes the proof of Lemma \ref{l:212}.
\end{pf}
\begin{lem}\label{l:213}
Let $A$, $B$ and $C_j$ be the constants in $(E_v)$, that is,
$$
 \frac{\partial\psi_\e}{\partial v} = A P\d_{x,\l} +  B\frac{\partial P\d_{x,\l}}{\partial\l} + \sum_{j=1}^n C_j  \frac{\partial P\d_{x,\l}}{\partial x_j},
$$
where $\psi_\e$ is defined by \eqref{psi}. Then we have the following estimates
\begin{align}
& A= O\left(\e^{1-\sigma}\right),\quad B=O\left(\e^{1-1/(n-4)}\right),\quad C_j=O\left(\e^{1-\sigma + 1/(n-4)}\right),\label{e:271}\\
& \frac{\partial A}{\partial\l}= O\left(\e^{1-\sigma + 1/(n-4)}\right),\quad  \frac{\partial B}{\partial\l}= O\left(\e\right),\quad  \frac{\partial C_j}{\partial\l}= O\left(\e^{1-\sigma + 2/(n-4)}\right).\label{e:272}
\end{align}
\end{lem}
\begin{pf}
By Lemma \ref{l:28}, we see that $A$, $B$ and $C_j$ satisfy
\begin{align}
&A||P\d_{x,\l}||^2 + B \left( \frac{\partial P\d_{x,\l}}{\partial\l}, P\d_{x,\l}\right) + \sum_{j=1}^n C_j \left( \frac{\partial P\d_{x,\l}}{\partial x_j}, P\d_{x,\l}\right)= O\left(\e^{1-\sigma}\right), \label{e:273}\\
&A\left(P\d_{x,\l},\frac{\partial P\d_{x,\l}}{\partial\l}\right)  + B \biggl |\biggl | \frac{\partial P\d_{x,\l}}{\partial\l}\biggr |\biggr |^2 + \sum_{j=1}^n C_j \left( \frac{\partial P\d_{x,\l}}{\partial x_j}, \frac{\partial P\d_{x,\l}}{\partial\l}\right)= O\left(\e^{1+\frac{1}{n-4}}\right), \label{e:274}\\
&A\left(P\d_{x,\l},\frac{\partial P\d_{x,\l}}{\partial x_i}\right) + B \left( \frac{\partial P\d_{x,\l}}{\partial\l}, \frac{\partial P\d_{x,\l}}{\partial x_i}\right) + \sum_{j=1}^n C_j \left( \frac{\partial P\d_{x,\l}}{\partial x_j},\frac{\partial P\d_{x,\l}}{\partial x_i}\right)\notag\\
&\qquad\qquad \qquad \qquad \qquad = O\left(\e^{1-\sigma -\frac{1}{n-4}}\right), \quad\mbox{for}\quad 1\leq i\leq n.\label{e:275}
\end{align}
Solving \eqref{e:273},..., \eqref{e:275}, we obtain \eqref{e:271}. Differentiating  \eqref{e:273},..., \eqref{e:275} with respect to $\l$ and using Lemma \ref{l:212}, we obtain
\begin{align}
&\frac{\partial A}{\partial\l}||P\d_{x,\l}||^2 + \frac{\partial B}{\partial\l} \left( \frac{\partial P\d_{x,\l}}{\partial\l}, P\d_{x,\l}\right) + \sum_{j=1}^n \frac{\partial C_j}{\partial\l} \left( \frac{\partial P\d_{x,\l}}{\partial x_j}, P\d_{x,\l}\right)= O\left(\e^{1-\sigma + \frac{1}{n-4}}\right), \label{e:276}\\
&\frac{\partial A}{\partial\l}\left(P\d_{x,\l},\frac{\partial P\d_{x,\l}}{\partial\l}\right)  + \frac{\partial B}{\partial\l} \biggl |\biggl | \frac{\partial P\d_{x,\l}}{\partial\l}\biggr |\biggr |^2 + \sum_{j=1}^n \frac{\partial C_j}{\partial\l} \left( \frac{\partial P\d_{x,\l}}{\partial x_j}, \frac{\partial P\d_{x,\l}}{\partial\l}\right)= O\left(\e^{1 +\frac{2}{n-4}}\right), \label{e:277}\\
&\frac{\partial A}{\partial\l}\left(P\d_{x,\l},\frac{\partial P\d_{x,\l}}{\partial x_i}\right) +\frac{\partial B}{\partial\l} \left( \frac{\partial P\d_{x,\l}}{\partial\l}, \frac{\partial P\d_{x,\l}}{\partial x_i}\right) + \sum_{j=1}^n \frac{\partial C_j}{\partial\l} \left( \frac{\partial P\d_{x,\l}}{\partial x_j},\frac{\partial P\d_{x,\l}}{\partial x_i}\right)\notag\\
&\qquad\qquad \qquad \qquad \qquad = O\left(\e^{1-\sigma}\right), \quad\mbox{for}\quad 1\leq i\leq n.\label{e:278}
\end{align}
Solving \eqref{e:276},..., \eqref{e:278}, we get \eqref{e:272}.
\end{pf}
\section{Proof of Theorems \ref{t:11}, \ref{t:12} and \ref{t:16}}
First, let us introduce some notations. For two constants $\b$ and $L$ such that $L>\b >0$ we define a set
\be \label{e:31}
D_\e= \{x \in \O \cap \ov{B_{\e^\b}(x_0)}\mid d(x,\partial\O)\geq \e^L\}.
\ee
For constants $0<C_0 < C_1$ we set
$$
\l_{C_i}^\e (x) = C_i \left(\frac{H(x,x)}{\e}\right)^{1/(n-4)}\quad i=0,1
$$
and we define the following set
\be \label{e:32}
M_{\e}=\{ (x,\l)\mid x\in D_\e,\,\, \l \in [\l_{C_0}^\e(x), \l_{C_1}^\e(x)]\}.
\ee
Constants $\b$, $L$ and $C_i$ will be determined later.  We now consider the following minimization problem
\be\label{e:33}
\inf \{\psi_\e(x,\l, v_\e)\mid (x,\l)\in M_\e\},
\ee
where $v_\e$ is defined in Proposition \ref{p:21}.
It is obvious that for small fixed $\e>0$ problem \eqref{e:33} has a minimizer $(x_\e,\l_\e)$. In order to prove that $(x_\e, \l_\e,v_\e)$ is a critical point of $\psi_\e$, we only need to prove that $(x_\e,\l_\e)$ is an interior point of $M_\e$.

\begin{pfn}{\bf Theorem \ref{t:11}}
We prove that if $\e>0$ is small enough, the minimizer $(x_\e,\l_\e)$ of \eqref{e:33} is an interior point of $M_\e$. First we show that if $C_0$ and $C_1$ are suitably chosen, then
\be\label{e:34}
\l_\e \in (\l_{C_0}^\e(x_\e), \l_{C_1}^\e(x_\e)).
\ee
Using Proposition \ref{p:21} and the fact that $(x_\e, \l_\e)$ is a minimum point of \eqref{e:33}, we obtain
\be\label{e:35}
\psi_\e(x_\e, \l_\e, v_\e) \leq \psi_\e(x_\e,\l, 0)\quad\mbox{for all}\quad \l\in[\l_{C_0}^\e(x_\e), \l_{C_1}^\e(x_\e)].
\ee
As in the proof of Proposition \ref{p:21}, we obtain
\begin{align}\label{e:36}
\psi_\e(x_\e, \l_\e, v_\e)&= \psi_\e(x_\e,\l_\e, 0) + O\left(||v_\e||^2\right)\notag\\
&=\psi_\e(x_\e,\l_\e, 0) + O\left(\frac{1}{\l_\e^2} + \e^2 + \frac{1}{\left(\l_\e d_\e\right)^{n-4+2\th}}\right).
\end{align}
It follows from Proposition \ref{p:24} that
\begin{align}\label{e:37}
\frac{c_1H(x_\e,x_\e)}{S_n \l_\e^{n-4}}&+ \frac{n-4}{n}\e\left(\log \l_\e^{\frac{n-4}{2}} + \frac{c_3}{S_n}\right) + O\left(\frac{1}{\l_\e^2}+ \e^2 \log^2\l_\e + \frac{1}{(\l_\e d_\e)^{n-4+2\th}}+ \frac{\e \log\l_\e}{(\l_\e d_\e)^{n-4}} \right)\notag\\
&\leq \frac{c_1H(x_\e,x_\e)}{S_n \l^{n-4}}+ \frac{n-4}{n}\e\left(\log \l^{\frac{n-4}{2}} + \frac{c_3}{S_n}\right) + O\left(\frac{1}{\l^2}+ \e^2 \log^2\l + \frac{1}{(\l d_\e)^{n-4+2\th}}\right)\notag\\
& + O\left(  \frac{\e \log\l}{(\l d_\e)^{n-4}}\right).
\end{align}
Since $x_\e\in D_\e$, we get  $\e^L\leq d_\e \leq \e^\b$. If we choose $\b$ satisfying
\be \label{e:38}
\b > \max \{\frac{1}{2}-\frac{1}{n-4}, 0\},
\ee
then there exists a $\g >0$, such that
\begin{align}
\frac{1}{\l^2}& \leq C\left(\frac{\e}{H(x_\e,x_\e)}\right)^{\frac{2}{n-4}}= O\left(\e^{\frac{2}{n-4}}d_\e^2\right)= O\left(\e^{2\b + \frac{2}{n-4}}\right)= O\left(\e^{1+\g}\right),\label{e:39}\\
\e^2\log^2\l& = O\left(\e^2 \log \left(\frac{1}{\e^{\frac{1}{n-4}}d_\e}\right)\right)=
 O\left(\e^2 \log \left(\frac{1}{\e^{\frac{1}{n-4}+L}}\right)\right)= O\left(\e^{1+\g}\right),\label{e:310}\\
\frac{1}{\l d_\e}&\leq C\left(\frac{\e}{H(x_\e,x_\e)}\right)^{\frac{1}{n-4}}\frac{1}{d_\e}= O\left(\e^{1/(n-4)}\right).\label{e:311}
\end{align}
Consequently, we have
\be\label{e:312}
\frac{1}{(\l d_\e)^{n-4+2\th}}= O\left(\e^{1+\g}\right),\quad \frac{\e \log\l}{(\l d_\e)^{n-4}}= O\left(\e^{1+\g}\right).
\ee
Inserting \eqref{e:39},..., \eqref{e:312} into \eqref{e:37}, we obtain
\be\label{e:313}
 \frac{c_1H(x_\e,x_\e)}{S_n \l_\e^{n-4}}+ \frac{n-4}{n}\e\log \l_\e^{\frac{n-4}{2}}\leq \frac{c_1H(x_\e,x_\e)}{S_n \l^{n-4}}+ \frac{n-4}{n}\e\log \l^{\frac{n-4}{2}}+  O\left(\e^{1+\g}\right).
\ee
Let
\be\label{e:314}
\l_\e = t_\e \left(\frac{H(x_\e,x_\e)}{\e}\right)^{1/(n-4)},\quad \l = t \left(\frac{H(x_\e,x_\e)}{\e}\right)^{1/(n-4)}.
\ee
We then have from \eqref{e:313}
\be\label{e:315}
 \frac{c_1}{S_n t_\e^{n-4}}+ \frac{n-4}{n}\log t_\e^{\frac{n-4}{2}}\leq \frac{c_1}{S_n t^{n-4}}+ \frac{n-4}{n}\log t^{\frac{n-4}{2}}+  O\left(\e^{\g}\right).
\ee
Since $t\mapsto\frac{c_1}{S_n t^{n-4}}+ \frac{n-4}{n}\log t^{\frac{n-4}{2}}$, $t>0$, attains its global minimum at
\be\label{e:316}
t^*=\left(\frac{2nc_1}{(n-4)S_n}\right)^{1/(n-4)},
\ee
we deduce from \eqref{e:315} that as $\e\to 0$, $t_\e \to t^*$. If we choose
$$
C_0= \frac{1}{2}\left(\frac{2nc
_1}{(n-4)S_n}\right)^{1/(n-4)},\quad C_1= \frac{3}{2}\left(\frac{2nc_1}{(n-4)S_n}\right)^{1/(n-4)},
$$
then, for $\e>0$ small, we obtain \eqref{e:34}. Now, it remains to prove that $x_\e$ is an interior point of $D_\e$. To this aim, let $\nu$ be the inward unit normal of $\partial\O$ at $x_0$. Let $z_\e= x_0 + \e\nu$ and fix $\l_\e^* \in (\l_{C_0}^\e(z_\e), \l_{C_1}^\e(z_\e))$. Since $d(z_\e,\partial\O)=\e$, we have
$$
\l_\e^* \sim \left(\frac{1}{\e\e^{n-4}}\right)^{1/(n-4)}=\e^{-(n-3)/(n-4)}.
$$
We have
\be\label{e:317}
\psi_\e(x_\e,\l_\e,v_\e)\leq \psi_\e(z_\e,\l_\e^*,0).
\ee
In view of Proposition \ref{p:24}, we have
\begin{align}\label{e:318}
\psi_\e(z_\e,\l_\e^*,0)&= \frac{S_n^{(p-1-\e)/(p+1-\e)}}{\left(K(z_\e)\right)^{2/(p+1-\e)}}\left[ 1+ O\left(\e\log(1/\e)\right)\right]\notag\\
&= \frac{S_n^{(p-1-\e)/(p+1-\e)}}{\left(K(x_0)\right)^{2/(p+1-\e)}}\left(1+O(\e)\right)\left[ 1+ O\left(\e\log(1/\e)\right)\right]\notag\\
&= \frac{S_n^{(p-1-\e)/(p+1-\e)}}{\left(K(x_0)\right)^{2/(p+1-\e)}}\left[ 1+ O\left(\e\log(1/\e)\right)\right].
\end{align}
Using \eqref{e:318} and Proposition \ref{p:24}, \eqref{e:317} becomes
\begin{align}\label{e:319}
\frac{S_n^{\frac{p-1-\e}{p+1-\e}}}{K(x_\e)^{\frac{2}{p+1-\e}}}&\left[ 1+ \frac{c_1H(x_\e,x_\e)}{S_n\l_\e^{n-4}}+ \frac{n-4}{n}\e \left(\log\l_\e^{\frac{n-4}{2}} + \frac{c_3}{S_n}\right) + O\left(\e^{1+\g}\right)\right]\notag\\
&\leq \frac{S_n^{(p-1-\e)/(p+1-\e)}}{\left(K(x_0)\right)^{2/(p+1-\e)}}\left[ 1+ O\left(\e\log(1/\e)\right)\right].
\end{align}
We now consider two steps.\\
{\it Step 1}. We claim that $x_\e \not\in \{x \mid
d(x_\e,\partial\O)= \e^L\}$. Arguing by contradiction, suppose
that $ d(x,\partial\O)= \e^L$. Then \be\label{e:320} \l_\e \geq
\l_{C_0}^\e(x_\e)\geq c\left(\frac{1}{d(x_\e,\partial\O)^{n-4}
\e}\right)^{1/(n-4)}\geq \frac{C}{\e^L}. \ee Since $K(x_\e)\leq
K(x_0) $ and using \eqref{e:320}, \eqref{e:319} implies

$$
\frac{(n-4)^2}{2n}L\e \log(1/\e) + O(\e)\leq C\e \log(1/\e),
$$
where $C$ is a positive constant independent of $L$. So we get a contradiction if $L$ is chosen large enough.\\
{\it Step 2}. We claim that $x_\e \not\in \partial B_{\e^\b}(x_0)$. Again arguing by contradiction we assume that  $x_\e \in \partial B_{\e^\b}(x_0)$. Then by assumption on $K$, we have
$$
\frac{1}{K(x_\e)^{2/(p+1-\e)}}\geq \frac{1}{\left(K(x_0)-a\e^{\b(2+\a)}\right)^{2/(p+1-\e)}}\geq \frac{1+a'\e^{\b(2+\a)}}{K(x_0)^{2/(p+1-\e)}},
$$
where $a'>0$. Hence, if we can choose $\b>0$ satisfying
\be\label{e:321}
\b(2+\a) <1,
\ee
then, using \eqref{e:319}, we obtain
\be\label{e:322}
a'\e^{\b(2+\a)}\leq C\e\log (1/\e),
\ee
which is impossible. Thus it remains to prove that we can choose a $\b>0$, such that \eqref{e:38} and \eqref{e:321} hold. We distinguish two cases: $(i) n\geq 7$ and $(ii) n=5, 6$. In the case $(i)$ since $\a\in [0,\frac{4}{n-6})$, we can choose $\b\in (1/2 - 1/(n-4), 1/2)$ satisfying $\b(2+\a)<1$. Finally, if $n=5,6$, we can take $\b>0$ sufficiently small such that \eqref{e:321} holds.
From Steps $1$ and $2$ we deduce that $x_\e$ is an interior point of $D_\e$.\\
 By construction, the corresponding $u_\e=P\d_{x_\e,\l_\e}+ v_\e$ is a critical point of $J_\e$, that is, $w_\e = (l_\e(u_\e)^{\frac{1}{p-1-\e}} u_\e$ satisfies
 \begin{eqnarray}\label{e:ew}
 \D ^2 w_\e =K |w_\e|^{\frac{8}{n-4} - \e}w_\e \mbox{ in }\O, \quad w_\e=\D
 w_\e=0 \mbox{ on } \partial \O
 \end{eqnarray}
with $|w_\e^-|_{L^{2n/(n-4)}(\O)}$ very small, where $w_\e
^-=\max(0,-w_\e)$. As in Proposition 4.1 of \cite{BEH}, we can prove
that $w_\e^-=0$. Thus, since $w_\e$ is a non-negative function
which satisfies \eqref{e:ew}, the strong maximum principle ensures
that $w_\e > 0$ on $\O$ and then $u_\e$ is a solution of
$(P_{\e})$. This ends the proof of our Theorem.
\end{pfn}

\begin{pfn}{\bf Theorem \ref{t:12}}
Since the proof of Theorem \ref{t:12} is similar to that of
Theorem \ref{t:11}, we only point out the necessary changes in the
proof.
Let $\d >0$ such that
$\forall x \in B_\d(x_0), K(x)\leq K(x_0) $. We consider the minimization problem. \be\label{e:323}
\inf\{\psi_\e(x,\l,v_\e)\mid x\in \ov{B_\d(x_0)}, \l\in [\e^{-\b},
\e^{-L}]\}, \ee in place of \eqref{e:33}, where $0< \b < L$ are
some constants to be determined later. Let $(x_\e,\l_\e)$ be
a minimizer of problem \eqref{e:323}. From
$\psi_\e(x_\e,\l_\e,v_\e)\leq \psi_\e(x_0,\l_\e, 0)$, and the fact that $x_0$ is a strict local maximum, we easily
derive that $x_\e \to x_0$. Next, we show that $L$ and $\b$ can be
chosen so that $\e^{-\b} < \l_\e < \e^{-L}$. On one hand, it
follows from Proposition \ref{p:24} that \be\label{e:324}
\psi_\e(x_0,\e^{-4},0)=
\frac{S_n^{\frac{p-1-\e}{p+1-\e}}}{K(x_0)^{\frac{2}{p+1-\e}}}
\left[ 1+\frac{(n-4)^2}{2n}\e\left(4\log(1/\e)
+\frac{2c_3}{(n-4)S_n}\right) + O\left(\e^{1+\s}\right)\right].
\ee On the other hand, by Lemma \ref{l:23} and \eqref{e:216}, we have
\begin{align*}
\int_\O &K(y)|P\d_{x,\l} +v_\e|^{p+1-\e}= \int_{B_\d(x_0)} K(y)|P\d_{x,\l} +v_\e|^{p+1-\e} + O\left(\l^{-n}+ ||v_\e||^{p+1-\e}\right) \\
&\leq  K(x_0)\int_{B_\d(x_0)} |P\d_{x,\l} +v_\e|^{p+1-\e} + O\left(\l^{-n}+ ||v_\e||^{p+1-\e}\right) \\
&=  K(x_0)\biggl [\int_\O |P\d_{x,\l}|^{p+1-\e}+ (p+1-\e)\int_\O |P\d_{x,\l}|^{p-\e}v_\e \\
&+ \frac{(p+1-\e)(p-\e)}{2}\int_\O |P\d_{x,\l}|^{p-1-\e}v_\e^2\biggr] + O\left(\l^{-n}+ ||v_\e||^{\min(3,p+1-\e)}\right) \\
&=K(x_0)\biggl[S_n -(p+1-\e)\frac{c_1H(x,x)}{\l^{n-4}}-\e S_n\left(\log\l^{\frac{n-4}{2}} + \frac{c_3}{S_n}\right)+ O\left(\e + \frac{1}{\l^{\frac{n-4}{2}+\th}}\right)||v_\e||\\
&+\frac{(p+1-\e )(p-\e )}{2}\int_\O P\d_{x,\l}^{p-1-\e }v_\e
^2\biggr] + O\left(\frac{\e \log\l }{\l ^{n-4}} + \frac{1}{\l
^{n-2}}+ \e^2\log^2\l + ||v_\e||^{\min(3,p+1-\e)}\right) .
\end{align*}
Clearly, the above estimate implies
\begin{align}\label{e:325}
\psi_\e(x,\l,v_\e)&\geq \frac{S_n^{\frac{p-1-\e}{p+1-\e}}}{K(x_0)^{\frac{2}{p+1-\e}}} \biggl[ 1+ \frac{c_1H(x,x)}{S_n\l^{n-4}} + \frac{(n-4)}{n}\e\left(\log\l^{\frac{n-4}{2}} +\frac{c_3}{S_n}\right) \notag \\
&+\rho ||v_\e||^2 + O\left(\e^2\log^2\l + \frac{1}{\l^{n-4+2\th}}+ \frac{\e\log\l}{\l^{n-4}} + \frac{1}{\l^{n-2}}\right)\biggr] \notag\\
&\geq
\frac{S_n^{\frac{p-1-\e}{p+1-\e}}}{K(x_0)^{\frac{2}{p+1-\e}}}
\left[1+\frac{c_1H(x,x)}{2S_n\l^{n-4}}+\frac{(n-4)}{2n}\e\left(\log\l^{\frac{n-4}{2}}
+\frac{c_3}{S_n}\right) \right].
\end{align}
Using the inequality $\psi_\e(x_\e,\l_\e,v_\e)\leq
\psi_\e(x_0,\e^{-4}, 0)$, we deduce from \eqref{e:324} and
\eqref{e:325} that \be\label{e:326}
\frac{c_1H(x_\e,x_\e)}{2S_n\l_\e^{n-4}}+\frac{n-4}{2n}\e
\log\l_\e^{\frac{n-4}{2}}\leq \frac{(n-4)^2}{2n}4\e \log(1/\e) +
O\left(\e^{1+\s}\right). \ee
As in the proof of Theorem \ref{t:11} we proceed in two steps.\\
{\it Step 1}. We claim that $ \l_\e < \e^{-L}$ for $L>0$ sufficiently large. Arguing by contradiction, suppose that $\l_\e= \e^{-L}$. Then it follows from \eqref{e:326} that
$$
\e L\log(1/\e) \leq 8\e\log(1/\e) + O\left(\e^{1+\s}\right),
$$
which is impossible if $L$ is large enough.\\
{\it Step 2}. $\l_\e= \e^{-\b}$ is impossible if $\b>0$ is small enough. Assuming that $\l_\e=\e^{-\b}$, we deduce from \eqref{e:326} that
$$
\e^{(n-4)\b}\leq C \e\log(1/\e),
$$
which is impossible if $\b$ is small enough and therefore our result follows.
\end{pfn}

\begin{pfn}{\bf Theorem \ref{t:16}}
Arguing by contradiction, suppose that $(P_\e)$ has  a solution of the form \eqref{e:18} and satisfying \eqref{e:19}. We start by showing that $\e\log\l_\e \to 0$ as $\e\to 0$. Indeed, multiplying $(P_\e)$ by $P\d_{x_\e,\l_\e}$ and integrating over $\O$, we obtain
\begin{align*}
\a_\e ||P\d_{x_\e,\l_\e}||^2&= \int_\O K(y)|\a_\e P\d_{x_\e,\l_\e} + v_\e|^{p-\e} P\d_{x_\e,\l_\e}\\
&=\a_\e^{p-\e}\int_\O K(y)P\d_{x_\e,\l_\e}^{p+1-\e} + O\left(||v_\e||\right)\notag.
\end{align*}
As in \eqref{e:213}, we have
\begin{align}\label{e:326'}
\int_\O K(y)P\d_{x_\e,\l_\e}^{p+1-\e} = \frac{ S_n K( x_\e )}{\l_\e ^{\e(n-4)/2}} (1+o(1)).
\end{align}
Consequently by \eqref{e:212}, \eqref{e:326'}  we have
\be\label{e:327}
\a_\e S_n= \frac{ \a_\e^{p-\e} S_n K( x_\e )}{\l_\e ^{\e(n-4)/2}} \left(1+o(1)\right)+ o(1),
\ee
where $o(1)\to 0$ as $\e\to 0$. Since $\a_\e\to K(x_0)^{(4-n)/8}$ and $x_\e\to x_0$ as $\e \to 0$, we deduce from \eqref{e:327} that $\e\log\l_\e \to 0$ as $\e\to 0$.

Next, we estimate $v_\e$. Multiplying $(P_\e)$ by $v_\e$ and integrating over $\O$, we obtain
\begin{align}\label{e:328}
||v_\e||^2&= \int_\O K(y)|\a_\e P\d_{x_\e,\l_\e} + v_\e|^{p-\e} v_\e\notag\\
&=\a_\e^{p-\e}\int_\O K(y)P\d_{x_\e,\l_\e}^{p-\e}v_\e +
(p-\e)\a_\e^{p-1-\e}\int_\O K(y)P\d_{x_\e,\l_\e}^{p-1-\e}v_\e^2 +
O\left(||v_\e||^{\min(3,p+1-\e)}\right)\notag\\
&=\a_\e^{p-\e}\int_\O K(y)P\d_{x_\e,\l_\e}^{p-\e}v_\e +
(p-\e)\a_\e^{p-1-\e}K(x_\e)\int_\O \d_{x_\e,\l_\e}^{p-1}v_\e^2 +
o\left(||v_\e||^2\right).
\end{align}
It follows from Proposition 3.4 in \cite{BE1} that there exists a $\rho >0$, such that
\begin{align}\label{e:329}
||v_\e||^2&-  (p-\e)\a_\e^{p-1-\e} K(x_\e)\int_\O \d_{x_\e,\l_\e}^{p-1}v_\e^2 \notag\\
&=||v_\e||^2-  p\int_\O \d_{x_\e,\l_\e}^{p-1}v_\e^2  +
o\left(||v_\e||^2\right)\geq \rho ||v_\e||^2.
\end{align}
Combining \eqref{e:328}, \eqref{e:329} and with the aid of Lemma
\ref{l:23} we get \be\label{e:330} ||v_\e|| =
O\left(\frac{|DK(x_\e)|}{\l_\e} + \e +\frac{1}{\l_\e^2} +
\frac{1}{\l_\e^{\th + (n-4)/2}}\right). \ee
We now assume that
$n\geq 7$. Multiplying $(P_\e)$ by $\partial
P\d_{x_\e,\l_\e}/\partial\l$ and integrating over $\O$, we derive
that
$$
\a_\e\left(P\d_{x_\e,\l_\e}, \frac{\partial P\d_{x_\e,\l_\e}}{\partial\l}\right) - \a_\e^{(p-\e)}
 \int_\O K(y)|P\d_{x_\e,\l_\e}+v_\e|^{p-\e} \frac{\partial P\d_{x_\e,\l_\e}}{\partial\l} = 0.
$$
Arguing as in the proof of Lemma \ref{l:26}, we easily arrive at
\be\label{e:331} \frac{c_2\D K(x_\e)}{n\l_\e^3} +
\frac{(n-4)S_nK(x_\e)\e}{2\l_\e} + O\left(\frac{1}{\l_\e^4}+
\frac{\e\log\l_\e}{\l_\e^3} + \frac{\e^2\log\l_\e}{\l_\e}+\frac{
|D K(x_\e)|^2}{\l_\e ^3}+\frac{\e^2}{\l_\e }\right)=0. \ee Since
$\D K(x_\e)>0$ and $\e\log\l_\e \to 0$ as $\e\to 0$, we get from
\eqref{e:331} that
$$
\frac{1}{\l_\e^3} + \frac{\e}{\l_\e} \leq 0,
$$
which is impossible.\\
Finally, we consider the case $n=6$. As in the case $n\geq 7$ we derive the following relation
$$
\left(\frac{ c_2\D K(x_\e) }{36K(x_\e)} -c_1
H(x_\e,x_\e)\right)\frac{1}{\l_\e^3} + \frac{S_n\e}{6\l_\e} +
o\left(\frac{\e}{\l_\e} + \frac{1}{\l_\e^3}\right) = 0,
$$
which contradicts the assumption $(ii)$. This completes the proof of Theorem \ref{t:16}.
\end{pfn}
\section{Proof of Theorems \ref{t:14}, \ref{t:15} and \ref{t:17}}
In this section, except in the proof of Theorems \ref{t:15} and \ref{t:17}, we always assume that $K$ satisfies the conditions in Theorem \ref{t:14}. We now  start by proving the following propositions.
\begin{pro}\label{p:41}
There exists an $\e_0>0$, such that for each $\e\in (0,\e_0]$, there exists a $C^1$-map $\l_\e : B_{\e^{1/L}} \to \R^+$, $x\mapsto \l_\e(x)$, such that $\l_\e(x)$ satisfies $(E_\l)$. Moreover, $\l_\e(x)= t_\e(x)\e^{-1/(n-4)}$ with
\be\label{e:41}
|t_\e(x)-t_0(x)| = O\left(\e^{\s}\right),
\ee
where $\s >0$ and 
$$
t_0(x)= \left(\frac{2nc_1H(x,x)}{(n-4)S_n}\right)^{1/(n-4)}.
$$
\end{pro}
\begin{pf}
Using Lemma \ref{l:26}, and the fact that \eqref{e:229} holds, we
obtain \be\label{e:42} \frac{\partial\psi_\e}{\partial
\l}(x,\l,v_\e) =
\frac{1}{\left(S_nK(x)\right)^{\frac{2}{p+1-\e}}}\left[
\frac{-c_1(n-4)H(x,x)}{\l^{n-3}}+ \frac{(n-4)^2S_n\e}{2n\l} +
O\left(\e^{1+\frac{1}{n-4}+\s}\right)\right]. \ee On the other
hand by Lemmas \ref{l:27},  \ref{l:213} we have \be\label{e:43}
 B\left( \frac{\partial ^2P\d_{x,\l}}{\partial\l^2}, v_\e\right) + \sum_{j=1}^n C_j \left( \frac{\partial ^2P\d_{x,\l}}{\partial x_j\partial \l}, v_\e\right)= O\left(\e^{1+\frac{1}{n-4}+ \frac{1-\s}{2}}\right).
\ee
Consequently, equation $(E_\l)$ is equivalent to
\be\label{e:44}
 -\frac{c_1(n-4)H(x,x)}{\l^{n-3}}+ \frac{(n-4)^2S_n\e}{2n\l} + O\left(\e^{1+\frac{1}{n-4}+\s}\right)=0.
\ee
Letting $\l_\e= t_\e \e^{-1/(n-4)}$, we deduce from \eqref{e:44} that
\be\label{e:45}
 -\frac{c_1H(x,x)}{t_\e^{n-3}}+ \frac{(n-4)S_n}{2n t_\e} + O\left(\e^{\s}\right)=0.
\ee
It is easy to see that \eqref{e:45} has a solution
$$
t_\e \in \left( \frac{1}{2}\left(\frac{2nc_1H(x,x)}{(n-4)S_n}\right)^{\frac{1}{n-4}}, \frac{3}{2}\left(\frac{2nc_1H(x,x)}{(n-4)S_n}\right)^{\frac{1}{n-4}}\right).
$$
This implies the existence of $\l_\e(x)$ satisfying $(E_\l)$. Next, we show that $\l_\e(x)$ is a $C^1$-map in $x$. To this aim, let
$$
F(\l)= \frac{\partial\psi_\e}{\partial \l} -B\left( \frac{\partial ^2P\d_{x,\l}}{\partial\l^2}, v_\e\right) - \sum_{j=1}^n C_j \left( \frac{\partial ^2P\d_{x,\l}}{\partial x_j\partial \l}, v_\e\right).
$$
Then it follows from Lemma \ref{l:212}, Lemma \ref{l:213} and Proposition \ref{p:211} that
\be\label{e:46}
F'(\l)= \frac{1}{(K(x)S_n)^{\frac{2}{p+1-\e}}}\left(\frac{(n-3)(n-4)c_1H(x,x)}{\l^{n-2}} - \frac{(n-4)^2S_n\e}{2n\l^2}\right) + O\left(\e^{1+\frac{2}{n-4}+ \s}\right)>0,
\ee
for all $\l \in \left( \frac{1}{2}\left(\frac{2nc_1H(x,x)}{(n-4)S_n\e}\right)^{\frac{1}{n-4}}, \frac{3}{2}\left(\frac{2nc_1H(x,x)}{(n-4)S_n\e}\right)^{\frac{1}{n-4}}\right)$. Consequently, the equation $(E_\l)$ has a unique solution in
$$
\l \in \left( \frac{1}{2}\left(\frac{2nc_1H(x,x)}{(n-4)S_n\e}\right)^{\frac{1}{n-4}}, \frac{3}{2}\left(\frac{2nc_1H(x,x)}{(n-4)S_n\e}\right)^{\frac{1}{n-4}}\right),
$$
and since all the terms in $(E_\l)$ are of $C^1$ with respect to $x$ and $\l$, we deduce that $\l_\e(x)$ is a $C^1$ map in $x$.

Now, let
$$
\Phi(t)= -\frac{(n-4)c_1H(x,x)}{t^{n-3}}+ \frac{(n-4)^2S_n}{2nt}.
$$
We then have
\be\label{e:47}
\Phi(t_0(x))=0,\quad \Phi(t_\e(x))= O\left(\e^\s\right).
\ee
Since $\Phi '(t_0(x)) >0$, it follows from \eqref{e:47} that $|t_\e(x)-t_0(x)|=  O\left(\e^\s\right)$ and this completes the proof of Proposition \ref{p:41}.
\end{pf}

Now, we consider the following maximization problem
\be\label{e:48} \sup \{\psi_\e(x,\l_\e(x), v_\e(x,\l_\e(x)))\mid
|x-x_0| \leq \e^{1/L}\}. \ee Then \eqref{e:48} has a maximizer
$x_\e\in\{x\mid |x-x_0|\leq \e^{1/L}\}$. In order to prove that
$x_\e$ is a critical point, we only need to prove that
$|x_\e-x_0|<\e^{1/L}$.
\begin{pro}\label{p:42}
Let $x_\e$ be a maximizer of \eqref{e:48}. Then there exists a $\s_2>0$, such that $|x_\e-x_0|^L=O\left(\e^{1+\s_2}\right)$. In particular, if $\e>0$ is small enough, $x_\e$ is an interior point of $B_{\e^{1/L}}(x_0)$.
\end{pro}
\begin{pf}
It follows from Lemma \ref{l:27}, Propositions \ref{p:24}, and \ref{p:41} that
\begin{align}\label{e:49}
\psi_\e(x,& \l_\e(x),v_\e(x,\l_\e(x)))\notag\\
&=  \frac{S_n^{\frac{p-1-\e}{p+1+\e}}}{K(x)^{\frac{2}{p+1-\e}}}\left[ 1+ \frac{c_1H(x,x)}{S_n\l_\e^{n-4}} + \frac{(n-4)\e}{n}\left(\log\l_\e^{\frac{n-4}{2}}+ \frac{c_3}{S_n}\right)+ O\left(\e^{1+\s}\right)\right].
\end{align}
Letting
$$
\l_\e(x):=t_\e(x) \e^{-1/(n-4)}=\left(t_0(x)+ O(\e^\s)\right)\e^{-1/(n-4)},$$
we deduce from \eqref{e:49} that
\begin{align}\label{e:410}
\psi_\e(x,\l_\e(x),v_\e(x,\l_\e(x)))
&= \frac{S_n^{\frac{p-1-\e}{p+1+\e}}}{K(x)^{\frac{2}{p+1-\e}}}\biggl[ 1+ \frac{c_1H(x,x)\e}{S_nt_0(x)^{n-4}} + \frac{(n-4)^2\e}{2n}\log t_0(x)\notag\\
&+ \frac{n-4}{2n}\e\log(1/\e)+ \frac{(n-4)c_3\e}{nS_n}+ O\left(\e^{1+\s}\right)\biggr].
\end{align}
Since $x_\e$ is a maximum of \eqref{e:48}, we have
$$
\psi_\e(x_\e,\l_\e(x_\e),v_\e(x_\e,\l_\e(x_\e)))\geq \psi_\e(x_0,\l_\e(x_0),v_\e(x_0,\l_\e(x_0))).
$$
This, together with \eqref{e:410} and the assumption
$$
K(x_\e)\geq K(x_0) + C_0|x_\e-x_0|^L,
$$
imply
\begin{align*}
|x_\e-x_0|^L & \leq C\e \left(\log H(x_\e,x_\e)- \log H(x_0,x_0)\right) + O\left(\e^{1+\s}\right)\notag\\
&=O\left(\e|x_\e-x_0|\right) + O\left(\e^{1+\s}\right).
\end{align*}
Hence $|x_\e-x_0|^L = O\left(\e^{1+\s_2}\right)$, where $\s_2$ is
a positive constant. Thus Proposition \ref{p:42} follows.
\end{pf}

\begin{pfn}{\bf Theorem \ref{t:14}}
We only need to prove that
$(x_\e,\l_\e(x_\e),v_\e(x_\e,\l_\e(x_\e)))$ satisfies $(E_x)$.
Indeed, we have by easy computations
\begin{align*}
0&=\frac{\partial\psi_\e}{\partial x_i}+ \frac{\partial\psi_\e}{\partial \l}\frac{\partial\l}{\partial x_i}+ \left(\frac{\partial\psi_\e}{\partial v}, \frac{\partial v}{\partial x_i}+ \frac{\partial v}{\partial\l}\frac{\partial\l}{\partial x_i}\right)\\
&=\frac{\partial\psi_\e}{\partial x_i}+\left[ B\left(\frac{\partial^2 P\d_{x,\l}}{\partial\l^2},v\right)+ \sum_{j=1}^n C_j \left(\frac{\partial^2 P\d_{x,\l}}{\partial\l \partial x_j},v\right)\right]\frac{\partial\l}{\partial x_i} +B\left(\frac{\partial P\d_{x,\l}}{\partial\l},\frac{\partial v}{\partial x_i}\right)\\
&+ \sum_{j=1}^n C_j \left(\frac{\partial P\d_{x,\l}}{ \partial x_j},\frac{\partial v}{\partial x_i}\right)+ \left[ B\left(\frac{\partial P\d_{x,\l}}{\partial\l},\frac{\partial v}{\partial\l}\right)+ \sum_{j=1}^n C_j \left(\frac{\partial P\d_{x,\l}}{ \partial x_j},\frac{\partial v}{\partial\l}\right)\right]\frac{\partial\l}{\partial x_i}\\
&=\frac{\partial\psi_\e}{\partial x_i}- B\left(\frac{\partial^2 P\d_{x,\l}}{\partial\l\partial x_i},v\right)- \sum_{j=1}^n C_j \left(\frac{\partial^2 P\d_{x,\l}}{\partial x_i \partial x_j},v\right).
\end{align*}
This obviously shows that $(E_x)$ holds and as in the proof of Theorem \ref{t:11}, we see that the corresponding $u_\e=P\d_{x_\e,\l_\e} + v_\e$ is a solution of $(P_\e)$.
\end{pfn}

\begin{pfn}{\bf Theorem \ref{t:15}}
 Theorem \ref{t:15} can be proved in exactly the same way as Theorem \ref{t:14}
\end{pfn}

 \begin{pfn}{\bf Theorem \ref{t:17}}
   Arguing by contradiction, let us suppose that $(Q_{\e})$ has  a solution of the form \eqref{e:18} and satisfying \eqref{e:19}. We start by showing that $\l_\e$ occurring in \eqref{e:18}
 satisfies $\l_\e^{\e(n-4)/2} \to 1$ as $\e\to 0$. Indeed, multiplying $(Q_{\e)}$ by $P\d_{x_\e,\l_\e}$ and integrating over $\O$, we obtain
\begin{align}
\a_\e ||P\d_{x_\e,\l_\e}||^2&= \int_\O K(y)|\a_\e P\d_{x_\e,\l_\e} + v_\e|^{p+\e} P\d_{x_\e,\l_\e}\notag\\
&=\a_\e^{p+\e}\int_\O K(y)P\d_{x_\e,\l_\e}^{p+1+\e} + O\left(\int_\O \d_{x_\e,\l_\e}^{p+\e}|v_\e|+\int_\O \d_{x_\e,\l_\e}|v_\e|^{p+\e} \right)\notag\\
&=\a_\e^{p+\e}\int_\O K(y)P\d_{x_\e,\l_\e}^{p+1+\e} + O\left(\l_\e^{\e(n-4)/2}\int_\O \d_{x_\e,\l_\e}^p|v_\e|+\l_\e^{\e(n-4)/2}\int_\O \d_{x_\e,\l_\e}^{1-\e}|v_\e|^{p+\e} \right)\notag\\
&=\a_\e^{p+\e}\int_\O K(y)P\d_{x_\e,\l_\e}^{p+1+\e} +
O\left(\l_\e^{\e(n-4)/2}||v||+\l_\e^{\e(n-4)/2}||v||^{p+\e}\right)\notag.
\end{align}
As in \eqref{e:213}, we have
\begin{align}\label{e:411}
 \int_\O K(y)P\d_{x_\e,\l_\e}^{p+1+\e} &= \int_\O K(y)\left(\d_{x_\e,\l_\e}-\varphi_{x_\e,\l_\e}\right)^{p+1+\e}\notag\\
&=\int_\O K(y)\d_{x_\e,\l_\e}^{p+1+\e}+O\left(\int_\O \d_{x_\e,\l_\e}^{p+\e}\varphi_{x_\e,\l_\e}\right)\notag\\
&=K(x_\e)c_0^\e \l^{\e(n-4)/2}\int_{\R^n}\d_{0,1}^{p+1+\e}+o(\l^{\e(n-4)/2})\notag\\
&=S_n K( x_\e )\l_\e ^{\e(n-4)/2} (1+o(1)).
\end{align}
Consequently by \eqref{e:212} and \eqref{e:411},  we have
  \be\label{e:412}
  \a_\e S_n=  \a_\e^{p+\e} S_n K( x_\e )\l_\e ^{\e(n-4)/2} \left(1+o(1)\right)+o(1).
  \ee
   Since $\a_\e\to K(x_0)^{(4-n)/8}$ and $x_\e\to x_0$ as $\e \to 0$, we deduce from \eqref{e:412} that $\l_\e^{\e(n-4)/2} \to 1$ as $\e\to 0$.

Next, we are going to estimate the $v_\e$-part of $u_\e$. Multiplying $(Q_{\e)}$ by $v_\e$ and
integrating over $\O$, we obtain
\begin{align}\label{e:413}
||v_\e||^2&= \int_\O K(y)|\a_\e P\d_{x_\e,\l_\e} + v_\e|^{p+\e} v_\e\notag\\
&=\a_\e^{p+\e}\int_\O K(y)P\d_{x_\e,\l_\e}^{p+\e}v_\e +
(p+\e)\a_\e^{p-1+\e}\int_\O K(y)P\d_{x_\e,\l_\e}^{p-1+\e}v_\e^2 \notag\\
&+ O\left(||v_\e||^3+\int_\O |v_\e|^{p+1+\e}\right).
\end{align}
 According to Lemmas 4.4 and 4.5 of \cite{BE}, we have
$$
\int_\O |v_\e|^{p+1+\e}=o(1)\quad \mbox{and}\quad |v|_{L^\infty(\O)}^{\e}=O(1),
$$
 therefore
  \begin{align}
  ||v_\e||^2-(p+\e)\a_\e^{p-1+\e}\int_\O K(y)P\d_{x_\e,\l_\e}^{p-1+\e}v_\e^2=\a_\e^{p+\e}\int_\O K(y)P\d_{x_\e,\l_\e}^{p+\e}v_\e +  + O(||v_\e||^{inf(3,
p+1)}).\notag
  \end{align}
 Observe that
\begin{align}
||v_\e||^2&-(p+\e)\a_\e^{p-1+\e}\int_\O K(y)P\d_{x_\e,\l_\e}^{p-1+\e}v_\e^2=||v_\e||^2-(p+\e)\a_\e^{p-1+\e}K(x_\e)\int_\O \d_{x_\e,\l_\e}^{p-1+\e}v_\e^2 +o(||v||^2)\notag\\
&=||v_\e||^2-(p+\e)\a_\e^{p-1+\e}K(x_\e)c_0^\e\l^{\e(n-4)/2}\int_\O
\d_{x_\e,\l_\e}^{p-1}v_\e^2+o(||v||^2).
\end{align}
Since $x_\e\to x_0$, $\a_\e\to K(x_0)^{(4-n)/8}$ and
$\l^{\e(n-4)/2}\to 1$ as $\e \to 0$, it follows from Proposition 3.4 of \cite{BE1} that
 \begin{align}
 ||v_\e||^2-(p+\e)\a_\e^{p-1+\e} K(x_\e)c_0^\e\l^{\e(n-4)/2}\int_\O \d_{x_\e,\l_\e}^{p-1}v_\e^2\geq\rho||v_\e||^2,
\end{align}
where $\rho$ is a positive constant independent of $\e$.\\
 As in Lemma \ref{l:23}, we have
\begin{align}\label{e:416}
\a_\e ^{p+\e }\int_\O K(y)P\d_{x_\e ,\l_\e }^{p+\e }v_\e
=O\left(\frac{|DK(x_\e )|}{\l_\e} + \e +\frac{1}{\l_\e ^2} +
\frac{1}{\l_\e ^{\th + (n-4)/2}}\right)
\end{align}
 Then we deduce from  \eqref{e:413} and \eqref{e:416} that
\begin{align}
  ||v_\e|| =O\left(\frac{|DK(x_\e)|}{\l_\e} + \e +\frac{1}{\l_\e^2} + \frac{1}{\l_\e^{\th +(n-4)/2}}\right).
\end{align}
Now, multiplying $(Q_{\e})$ by $\partial P\d_{x_\e,\l_\e}/\partial\l$
and integrating over $\O$, we derive that
$$
\a_\e\left(P\d_{x_\e,\l_\e}, \frac{\partial
P\d_{x_\e,\l_\e}}{\partial\l}\right) -
 \int_\O K(y)|\a_\e P\d_{x_\e,\l_\e}+v_\e|^{p+\e} \frac{\partial P\d_{x_\e,\l_\e}}{\partial\l} = 0.
$$
Arguing as in the proof of Lemma \ref{l:26}, we easily arrive at
\begin{align}\label{e:418}
 -\frac{c_2\D K(x_\e)}{n^2 K(x_\e)\l_\e^3} + \frac{ c_1H(x_\e,x_\e)}{\l_\e^{n-3}}
+\frac{(n-4)S_n\e}{2n\l_\e} + O\left(\frac{1}{\l_\e^4}+
\frac{\e\log\l_\e}{\l_\e^3} + \frac{\e^2\log\l_\e}{\l_\e}\right)\notag \\
+\left(\frac{ |D K(x_\e)|^2}{\l_\e ^3}+\frac{\e^2}{\l
}+\frac{1}{\l_\e^{n-3+2\th}}+\frac{\e \log
\l_\e}{\l_\e^{n-3}}\right)=0.
\end{align}
For $n=5$, it follows from \eqref{e:418} that
  $$
  \frac{c_1 H(x_\e,x_\e)}{\l_\e^2}+\frac{(n-4)S_n\e}{2n\l_\e} + o\left(\frac{1}{\l_\e^2}+\frac{\e}{\l_\e }\right)=0, 
$$
which is impossible.\\
For $n=6$, we derive form \eqref{e:418} the following relation 
$$
\left(-\frac{c_2\D K(x_\e)}{36K(x_\e)} + c_1 H(x_\e,x_\e)\right)
\frac{1}{\l_\e^3}+\frac{S_n\e}{6\l_\e} + o\left(\frac{1}{\l_\e^3}+
\frac{\e}{\l_\e}\right)=0, 
$$
 which is a contradiction with the assumption $(ii)$.\\
Finally, for $n\geq 7$,  we derive the following
relation 
$$
 -\frac{c_2\D K(x_\e)}{n^2K(x_\e)\l_\e^3}
+\frac{(n-4)S_n\e}{2n\l_\e} + o\left(\frac{1}{\l_\e^3}
+\frac{\e}{\l_\e }\right)=0,
 $$
 which contradicts the assumption
$(iii)$. This completes the proof of Theorem \ref{t:17}.
\end{pfn}
\section{Appendix}
In this appendix, we collect the integral estimates which are needed in Section 2.
\begin{lem}\label{l:22}
Suppose that $\l d(x,\partial\O)\to + \infty$, $\e \log \l
\to 0$ as
$\e\to 0$.\\ Then the following estimates hold 
\begin{align*}
1.\quad \int_\O K(y) \d_{x,\l}^{p+1-\e}(y) dy &= K(x) S_n + \frac{c_2\D K(x)}{2n\l ^2}  - \e K(x)S_n \left(\log \l^{\frac{n-4}{2}} + \frac{c_3}{S_n}\right)\quad\quad\quad\qquad\\
&+ O\left( \sum_{j=3}^{n-4} \frac{|D^jK(x)|}{\l^j} +
\frac{1}{\l^{n-3}}+\frac{\e \log \l}{\l^2} + \e^2 \log^2 \l +
\frac{1}{(\l d)^n}\right),
\end{align*}

\quad\quad where $S_n=\int_{\R^n}\d_{o,1}^{p+1}$, $c_2=
\int_{\R^n}|y|^2 \d_{o,1}^{p+1} dy$, $c_3=
\int_{\R^n}\d_{o,1}^{p+1}\log \d_{o,1}(y)dy$ and $d= d(x,
\partial\O )$.
$$
2.\quad \int_\O K(y) \d_{x,\l}^{p-\e}\var_{x,\l} =
\frac{c_1K(x)H(x,x)}{\l^{n-4}} + O\left(\frac{\e \log \l}{(\l
d)^{n-4}}+ \frac{1}{(\l d)^{n-2}}\right),\qquad\quad\quad\quad\quad\quad\qquad
$$
\qquad  where $c_1= c_0^{\frac{2n}{n-4}}\int_{\R^n} \frac{dy}{(1+|y|^2)^{(n+4)/2}}$ and $c_0$ is defined in \eqref{e:11}.
\begin{align*}
3.\quad \int_\O K(y) \d_{x,\l}^{p-\e}(y)\frac{\partial\d_{x,\l}}{\partial\l} dy &=- \frac{K(x)(n-4)^2 S_n\e}{4n\l} - \frac{(n-4)}{2n^2}c_2\frac{\D K(x)}{\l^3}\quad\qquad\qquad\qquad\,\,\\
&+ O\left( \frac{\e^2 \log \l}{\l}+\frac{\e \log
\l}{\l^3}+\sum_{j=3}^{n-4} \frac{|D^jK(x)|}{\l^{j+1}} +
\frac{1}{\l^{n-2}} + \frac{1}{\l(\l d)^n}\right).
\end{align*}
$$
4.\quad \int_\O K(y) \d_{x,\l}^{p-1-\e}\var_{x,\l}
\frac{\partial\d_{x,\l}}{\partial\l}=
-\frac{(n-4)^2}{2(n+4)}\frac{c_1K(x)H(x,x)}{ \l^{n-3}} +
O\left(\frac{\e \log \l}{\l (\l d)^{n-4}}+ \frac{1}{\l (\l
d)^{n-2}}\right).
$$
$$
5.\quad \int_\O K(y) \d_{x,\l}^{p-\e} \frac{\partial\var_{x,\l}}{\partial\l}= -\frac{(n-4)}{2}\frac{c_1K(x)H(x,x)}{ \l^{n-3}} + O\left(\frac{\e \log \l}{\l (\l d)^{n-4}}+ \frac{1}{\l (\l d)^{n-2}}\right).\qquad\quad
$$
\end{lem}
\begin{pf} 
Using the fact that $\d_{x,\l}^{-\e}= 1 -\e \log \d_{x,\l} + O \left(\e^2 \log^2 \l\right)$, we obtain
\begin{align*}
\int_\O K(y) \d_{x,\l}^{p+1-\e}(y) dy &= \int_\O K(y)\d_{x,\l}^{p+1}( 1- \e \log \d_{x,\l}) + O\left(\int_\O \d_{x,\l}^{p+1}(\e\log\l)^2\right)\\
&=\int_{B(x,d)} K(y) \d_{x,\l}^{p+1} -\e \log \l^{\frac{n-4}{2}}\int_{B(x,d)} K(y) \d_{x,\l}^{p+1}\\
&-\e \int_{B(x,d)} K(y) \d_{x,\l}^{p+1}\log \left(\frac{c_0}{(1+\l^2|y-x|^2)^{(n-4)/2}}\right)\\
& + O\left(\e^2 \log^2 \l + (\l d)^{-n}\right).
\end{align*}
Thus, using Taylor's expansion, we easily derive Claim 1.\\
Now, using (see \cite{BH})
$$
\var_{x,\l}= c_0 \frac{H(x,y)}{\l^{(n-4)/2}} + O\left( \frac{1}{\l^{n/2}d^{n-2}}\right),
$$
we derive that
\begin{align*}
\int_\O K(y)\d_{x,\l}^{p-\e}\var_{x,\l}&= \int_\O \d_{x,\l}^{p-\e}\frac{c_0K(y)H(x,y)}{\l^{(n-4)/2}} + O\left(\frac{1}{(\l d)^{n-2}}\right)\\
&=\int_{B(x,d)} \d_{x,\l}^{p}\frac{c_0 K(y)H(x,y)}{\l^{(n-4)/2}} + O\left(\frac{\e \log \l}{(\l d)^{n-4}} + \frac{1}{(\l d)^{n-2}}\right)\\
&=\frac{c_0 K(x)H(x,x)}{\l^{(n-4)/2}}\int_{\R^n} \d_{x,\l}^{p} +
O\left(\frac{\e \log \l}{(\l d)^{n-4}} + \frac{1}{(\l
d)^{n-2}}\right)
\end{align*}
and therefore Claim 2  follows.\\
To prove Claim 3, we use again Taylor's expansion and we thus obtain
\begin{align*}
\int_\O K(y)\d_{x,\l}^{p-\e}\frac{\partial\d_{x,\l}}{\partial\l}&= K(x)\int_{\R^n} \d_{x,\l}^{p-\e} \frac{\partial\d_{x,\l}}{\partial\l} + \frac{\D K(x)}{2n} \int_{\R^n} \d_{x,\l}^{p-\e}|y-x|^2 \frac{\partial\d_{x,\l}}{\partial\l}\\
&+ O\left(\sum_{j=3}^{n-4} \frac{|D^jK(x)|}{\l^{j+1}} + \frac{1}{\l^{n-2}} + \frac{1}{\l (\l d)^n}\right)\\
&= \frac{K(x)}{p+1-\e}\frac{\partial}{\partial\l}\left(\frac{1}{\l^{\e \frac{n-4}{2}}}\int_{\R^n}\d_{o,1}^{p+1-\e} \right)\\
& + \frac{\D K(x)}{2n(p+1-\e)}\frac{\partial}{\partial\l}\left(\frac{1}{\l^{2 +\e \frac{n-4}{2}}}\int_{\R^n}\d_{o,1}^{p+1-\e}|y|^2 \right)\\
&+  O\left(\sum_{j=3}^{n-4} \frac{|D^jK(x)|}{\l^{j+1}} + \frac{1}{\l^{n-2}} + \frac{1}{\l (\l d)^n}\right)\\
&= -K(x)\e \frac{(n-4)^2S_n}{4n \l}(1+O(\e))(1+O(\e \log\l ))
\end{align*}
\begin{align*}
\qquad\qquad 
& -c_2\frac{\D K(x)(n-4)}{2n^2\l^3}(1+O(\e))(1+O(\e \log \l))\\
&+  O\left(\sum_{j=3}^{n-4} \frac{|D^jK(x)|}{\l^{j+1}} + \frac{1}{\l^{n-2}} + \frac{1}{\l (\l d)^n}\right).
\end{align*}
Thus Claim 3  follows.\\
Now we are going to prove Claim 4. To this aim, we write
\begin{align*}
\int_\O K(y)\d_{x,\l}^{p-1-\e}\var_{x,\l}\frac{\partial\d_{x,\l}}{\partial\l}&=
\int_\O K(y)\d_{x,\l}^{p-1}\var_{x,\l}\frac{\partial\d_{x,\l}}{\partial\l} + O\left(\frac{\e \log \l}{\l} \int_\O \d_{x,\l}^{p}\var_{x,\l}\right)\\
&=
c_0\frac{K(x)H(x,x)}{\l^{\frac{n-4}{2}}}\int_{B(x,d)}\d_{x,\l}^{p-1}\frac{\partial\d_{x,\l}}{\partial\l}
+ O\left(\frac{1}{\l(\l d)^{n-2}}+ \frac{\e \log \l}{\l (\l
d)^{n-4}}\right)
\end{align*}
and we can thus easily derive Claim 4.\\
Lastly, we have
\begin{align*}
\int_\O K(y)\d_{x,\l}^{p-\e}\frac{\partial\var_{x,\l}}{\partial\l}&=
\int_\O K(y)\d_{x,\l}^{p}\frac{\partial\var_{x,\l}}{\partial\l} + O\left(\frac{\e \log \l}{\l (\l d^2)^{\frac{n-4}{2}}} \int_\O \d_{x,\l}^{p}\right)\\
&=- c_0\frac{n-4}{2}\int_\O \frac{K(y)H(x,y)}{\l^{\frac{n-2}{2}}}\d_{x,\l}^{p}  + O\left(\frac{1}{\l^{\frac{n+2}{2}}d^{n-2}}\int_\O \d_{x,\l}^p + \frac{\e \log \l}{\l (\l d)^{n-4}}\right)
\end{align*}
and thus Claim 5 follows. The proof of Lemma \ref{l:22} is thereby completed.
\end{pf}
\begin{lem}\label{l:23}
Let $k$ be the biggest positive integer satisfying $k\leq (n-4)/2$. Thus, for any $v\in E_{x,\l}$, we have
$$
1.\quad \int_\O K(y) P\d_{x,\l}^{p-\e}v = O\left(\e + \sum_{j=1}^{k} \frac{|D^jK(x)|}{\l^{j}} + \frac{1}{\l^{k+1}} + \frac{1}{(\l d)^{\frac{n-4}{2}+ \th}}\right)||v||,\qquad\qquad\quad
$$
$$
2.\quad \int_\O K(y) P\d_{x,\l}^{p-1-\e}v \frac{\partial
P\d_{x,\l}}{\partial\l}=  O\left(\frac{\e}{ \l} + \sum_{j=1}^k
\frac{|D^jK(x)|}{\l^{j+1}}+ \frac{1}{\l^{k+2}}+ \frac{1}{\l (\l
d)^{\frac{n-4}{2} + \th}}\right)||v||,
$$
$$
3.\quad \int_\O K(y) P\d_{x,\l}^{p-1-\e}v \frac{\partial
P\d_{x,\l}}{\partial x}=  O\left(\e\l + \sum_{j=1}^k
\frac{|D^jK(x)|}{\l^{j-1}}+ \frac{1}{\l^{k}}+ \frac{\l}{ (\l
d)^{\frac{n-4}{2} + \th}}\right)||v||,\quad
$$
where $\th$ is a positive constant.
\end{lem}
\begin{pf}
We observe that
\begin{align*}
\int_\O K(y) P\d_{x,\l}^{p-\e}v &=\int_\O K(y) \d_{x,\l}^{p-\e}v + O\left(\int_\O \d_{x,\l}^{p-1-\e}|\var_{x,\l}| |v|\right) \\
&=K(x)\int_\O \d_{x,\l}^p(1-\e \log (c_0\l^{(n-4)/2})v + O\left(\sum_{j=1}^k \frac{|D^jK(x)|}{\l^j} + \frac{1}{\l^{k+1}}\right)||v||\\
&+ O\left( \e \int_\O \d_{x,\l}^p \log ( 1+\l^2|x-a|^2)|v|\right) + O\left(\frac{||v||}{(\l d)^{\frac{n-4}{2}+\th}}\right)
\end{align*}
and thus Claim 1 follows.\\
We also observe that
\begin{align*}
\int_\O K(y) P\d_{x,\l}^{p-1-\e}&v \frac{\partial
P\d_{x,\l}}{\partial\l}= \int_\O K(y) P\d_{x,\l}^{p-1-\e}v
\frac{\partial \d_{x,\l}}{\partial\l} + O\left(\int_\O
\d_{x,\l}^{p-1-\e}|v| |\frac{\partial
\var_{x,\l}}{\partial\l}|\right)\\&= \int_\O K(y)
\d_{x,\l}^{p-1-\e}v \frac{\partial \d_{x,\l}}{\partial\l} +
O\left(\frac{1}{\l}\int_\O
\d_{x,\l}^{p-1-\e}|v|\var_{x,\l}\right) + O\left(\frac{||v||}{\l
(\l d)^{\frac{n-4}{2}+\th}}\right).
\end{align*}
Thus, using Taylor's expansion, we easily derive Claim 2.\\
In the same way, we can prove Claim 3 and therefore the proof of our lemma is completed.
\end{pf}

\end{document}